\documentclass{article}
\usepackage[utf8]{inputenc}
\usepackage[T1]{fontenc}
\usepackage{graphicx} % Required for inserting images
\usepackage{caption}
\usepackage{subcaption}

\usepackage[letterpaper, total={6in, 8in}]{geometry}
\usepackage{amsmath,amssymb,amscd,amsthm,amsfonts, bm}
\usepackage{hyperref}
\hypersetup{
    colorlinks=true,
    linkcolor=blue,
    filecolor=magenta,      
    urlcolor=cyan,
    pdftitle={Overleaf Example},
    pdfpagemode=FullScreen,
    }
\usepackage[section]{placeins}
\usepackage{parskip}
\usepackage[colorinlistoftodos,bordercolor=orange,backgroundcolor=orange!20,linecolor=orange,textsize=scriptsize]{todonotes}
\usepackage{amsthm}
\usepackage{verbatim}
\usepackage{thmtools}
\theoremstyle{plain}
\usepackage{pgf,tikz}
\usepackage{dirtytalk}
\usepackage{pdfpages}
\usepackage{mathrsfs}
\usepackage{enumitem}   
\usepackage{setspace}
\usepackage{imakeidx}
\usepackage{multirow}
\usepackage{layout}
\usepackage{lipsum}
\usepackage{thmtools}
\usepackage{thm-restate}
\usepackage{hyperref}

\usepackage{cleveref}

\usetikzlibrary{fit}

\makeindex

\usetikzlibrary{arrows}
\usepackage{authblk}
\usepackage[labelformat=simple]{subcaption}

\usepackage{
hyperref,%\autoref
}

\usepackage{textcomp}
\usepackage{soul}
\usepackage{color}
\usepackage{xcolor}

\makeatletter
\DeclareFontFamily{U}{tipa}{}
\DeclareFontShape{U}{tipa}{m}{n}{<->tipa10}{}
\newcommand{\arc@char}{{\usefont{U}{tipa}{m}{n}\symbol{62}}}%

\newcommand{\arc}[1]{\mathpalette\arc@arc{#1}}

\newcommand{\arc@arc}[2]{%
  \sbox0{$\m@th#1#2$}%
  \vbox{
    \hbox{\resizebox{\wd0}{\height}{\arc@char}}
    \nointerlineskip
    \box0
  }%
}
\makeatother

\definecolor{qqzzqq}{rgb}{0.,0.6,0.}
\definecolor{qqwuqq}{rgb}{0.,0.39215686274509803,0.}
\definecolor{ffqqqq}{rgb}{1.,0.,0.}
\definecolor{ududff}{rgb}{0.30196078431372547,0.30196078431372547,1.}

\newcommand{\Z}{\mathbb{Z}}

\newcommand{\bn}{\mathsf{\beta}}
\newcommand{\Bor}{\mathsf{DG}}
\newcommand{\Br}{\mathsf{Br}}

\newcommand\restr[2]{{% we make the whole thing an ordinary symbol
  \left.\kern-\nulldelimiterspace % automatically resize the bar with \right
  #1 % the function
  \littletaller % pretend it's a little taller at normal size
  \right|_{#2} % this is the delimiter
  }}

\newcommand{\littletaller}{\mathchoice{\vphantom{\big|}}{}{}{}}

\usepackage{subfiles}

\DeclareMathOperator{\verti}{vert}

\DeclareMathOperator{\diam}{\mathrm{diam}}
\newcommand{\St}{\widetilde S}

\newtheorem{theorem}{Theorem}
\newtheorem{definition}{Definition}
\newtheorem{conjecture}{Conjecture}
\newtheorem{question}{Question}

\newtheorem{lemma}{Lemma}
\newtheorem{proposition}{Proposition}

\newtheorem{remark}{Remark}

\title{On Diameter Graphs Problems:\\ A proof of Brass's Conjecture}
\author{Gyivan Lopez-Campos \\ 
\bigskip
\small\textit{Dedicated to Professor Jesús De Loera on his 60th birthday}}
\affil{University of California, Davis}
\date{\today}

\begin{document}

\maketitle
\begin{abstract}
The diameter graph of a finite set $S\subset\mathbb{R}^d$ has vertex set $S$ and two points are adjacent when their distance equals $\operatorname{diam} S$.
In this paper we study the following question: given $k$ and $d$, is there a constant
$c=c(k,d)$ such that from every finite point set $S\subset\mathbb{R}^d$ one can
remove at most $c$ points so that the subgraph of the diameter graph induced
by the remaining points is $k$-colorable? Brass conjectured that such a
constant exists for $d=3$ and $k=3$. We prove this conjecture by showing that at most one point needs to be removed. 

We also answer a question of
Dol'nikov, showing that for every positive integer $m$, there exist a finite point set
$S\subset\mathbb{R}^3$ in which every subset of smaller diameter has less than $(|S|-m)/2$ points. 

%In this paper we study some problems related with diameter graphs. We 
\end{abstract}

\section{Introduction}

Let $S\subset\mathbb{R}^d$ be a finite set of points in the $d$-dimensional real Euclidean space and let $\diam (S):=\sup_{x,y\in S}\mid\mid x-y\mid \mid$ be its Euclidean diameter. The \emph{diameter graph} of $S$ (also known as the \emph{largest-distance graph} or the \emph{Borsuk graph}) is the graph $\Bor (S)$ whose set of vertices is $S$, and two vertices are adjacent if their Euclidean distance equals $\diam (S)$. 

Diameter graphs play a central role in problems about diameters, such as the Borsuk partition problem, stated by Borsuk in \cite{borsuk1933drei}, the frequent largest distance problem proposed by Hopf and Pannwitz in \cite{hopf1934aufgabe} in dimension two and later posed by Vázsonyi in dimension three, and the Schur's Conjecture \cite{schur2003number}, proved in \cite{kupavskii2017proof}.

Recall that a graph is \emph{$k$-colorable} if all its vertices can be colored with $k$ colors so that any two adjacent vertices have different colors. In this paper, we address the following type of question:
\bigskip

\begin{question}\label{Question:Main}
             Let $k$ and $d$ be fixed positive integers. Is there a number of points, depending only on $k$ and $d$, that we can remove from any finite point set in $d$-dimensional Euclidean space, 
         such that the subgraph of the diameter graph induced by the remaining points is $k$-colorable?
\end{question}

This question is motivated by a problem posed by Brass in \cite{brassopen}, in which he asked about the case when $d=3$ and $k=2$. We believe this problem was proposed as attempt to address a more general question raised by Dol'nikov:
\bigskip

\begin{question} (Question 2, \cite{dol2000some}) \label{Question:Dolnikov}
    Can we choose such universal $m\in \mathbb{N}$, that for an arbitrary finite set $S$ in $\mathbb{R}^3$ there exist $U$ such that $U\subset S$, $|U|\geq (|S|-m)/2$, and $\diam (U)<\diam (S)$?
\end{question}

Notice that a positive answer to Brass’s problem would imply that, by removing at most a fixed number of points, one could obtain a bipartite diameter graph. Consequently, the answer to Question \ref{Question:Dolnikov} would also be positive.

However, Alon constructed in \cite{alon2003problems} a family of finite point sets that provides a negative answer to Brass’s problem, but leaving Question \ref{Question:Dolnikov} unanswered. Indeed, it is not difficult to show that the condition in Question~\ref{Question:Dolnikov} holds for this family with $m=4$.

Relaying on recent ideas introduced by Damási in  \cite{damasdi2026geometric1} and \cite{damasdi2026geometric}, in this paper we construct a new family of finite point sets that does give a negative answer to Question \ref{Question:Dolnikov}.
\bigskip

\begin{theorem}\label{Thm:Dolnikov}
    For every $m\in\mathbb{N}$, there exist a finite set $S\subset\mathbb{R}^3$, such that every $U\subset S$ with $\diam (U)<\diam (S)$ satisfies $|U|< (|S|-m)/2$.
\end{theorem}

Returning to Question \ref{Question:Main}, although no such constant exists for $k=2$ and $d=3$, Brass continued to ask questions this type. In the book \emph{Research Problems in Discrete Geometry} \cite[p.~231]{brass2005research}, he conjectured that, for $d=3$ and $k=3$, such a constant does exist. The conjecture is stated as follows.
\bigskip
\begin{conjecture}[Brass \cite{brass2005research}, 2005]
     From any finite point set in three-dimensional Euclidean space, one can remove $O(1)$ points with the property that the subgraph of the largest-distance graph induced by the remaining points is three-colorable.
\end{conjecture}

In this paper we prove that such a number does exist. More precisely, in Section \ref{Sec:BrassConjecture} we show the following theorem:
\bigskip
\begin{theorem}\label{Thm:Brass}
         From any finite point set in three-dimensional Euclidean space, one can remove at most one point with the property that the subgraph of the diameter graph induced by the remaining points is three-colorable.
\end{theorem}

The outline of the paper is as follows. In Section \ref{Sec:Preliminaries} we present the preliminaries. In Section \ref{Sec:Main} we prove Theorem \ref{Thm:Dolnikov} and Theorem \ref{Thm:Brass}. In Section \ref{Sec:Brasstype}, we address some cases of Question \ref{Question:Main}. Finally, in Section \ref{Sec:Final} we summarize what is known about Question \ref{Question:Main} and the open cases.

\textbf{Acknowledgments:} We would like to thank Professor Jesús Jeronimo-Castro for introducing us to the book  \emph{Research Problems in Discrete Geometry} and Dol'nivok's paper \cite{dol2000some}. This paper was partially supported by NSF grants DMS-2348578 and DMS-2434665.

\section{Preliminaries}\label{Sec:Preliminaries}

\subsection{Independence number, critical graphs and the Borsuk number}

In this paper, we assume familiarity with graph theory. The notation used is the same as that in \cite{bondy2008graph}. We refer the reader to this book for more details, proofs, and additional information on graph theory. However, to make this paper as self-contained as possible, we recall the notions that are crucial for understanding it. 

Let $G=(V,E)$ be a graph. A set $U\subset V$ is an \emph{independent set} of $G$ if no two vertices of $U$ are joint by an edge of $G$. The \emph{independence number} of $G$, denoted by $\alpha (G)$, is the maximum size of an independent set of $G$.

Notice that for a finite point set $S\subset\mathbb{R}^3$, a  subset $U\subset S$ is independent in $\Bor(S)$ if and only if $\diam (U)<\diam (S)$. Then, we can rephrase Theorem \ref{Thm:Dolnikov} as follows:

\bigskip
\begin{theorem}[Theorem \ref{Thm:Dolnikov} rephrased]\label{Thm:DolnikovRephrase}
        For every $m\in\mathbb{N}$, there exists a finite point set $S\subset\mathbb{R}^3$, such that $\alpha(\Bor(S))< (|S|-m)/2$.
\end{theorem}

Let $k$ be a positive integer. The \emph{chromatic number} $\chi (G)$ of the graph $G$ is the minimum $k$ such that $G$ is $k$-colorable. We say that the graph $G$ is \emph{$k$-critical} if $G$ has chromatic number $k$ and for all $v\in V$, the induced subgraph $G-\{v\}$ has chromatic number strictly less than $k$. A straightforward result on critical graphs is the following:\\
\smallskip

\begin{proposition}\label{Prop:critical}
    Every $k$-chromatic graph contains a $k$-critical subgraph.
\end{proposition}

The \emph{Borsuk number} $\bn(K)$ of a bounded set $K\subset \mathbb{E}^d$ is the smallest number of parts in which $K$ can be decomposed, such that each part has smaller diameter than $\diam (K)$. Clearly, if $K$ is finite, $\bn(K)=\chi(\Bor(K))$. 

The Borsuk number is invariant under translations,  rotations,  symmetries and uniform scalings.  Recall that a \emph{uniform scaling} is a linear transformation that enlarges (increases) or shrinks (diminishes) objects by a scale factor that is the same in all directions.  Therefore,  when we work on the Borsuk number we are just interested in the ``shape'' of the body and not in its embedding or  scale. 

By the previous observation, in this paper we assume without loss of generality that point sets have unit diameter unless otherwise stated.

The \emph{Borsuk number} of the dimension $d$ is defined as:
\[
\bn_d:=\sup\{\bn(K)|K\subset \mathbb{R}^d \text{bounded}\}.
\]

In 1933, Karol Borsuk asked in \cite{borsuk1933drei} whether $\bn_d=d+1$, raising a very interesting problem known as \emph{the Borsuk partition problem} or the \emph{Borsuk conjecture}. The Borsuk partition problem is known to have a positive answer for $d=2$ (see \cite{borsuk1933drei}) and $d=3$ (see \cite{perkal1947subdivision}, \cite{eggleston1955covering} and \cite{grunbaum1957simple}), and a negative answer for $d\geq 63$. 

In 1993, Kahn and Kalai were the first to find a counterexample \cite{kahn1993counterexample}, establishing that the answer is negative for $d> 2014$ and for $d=1325$. Since then, several mathematicians have found more counterexamples in progressively smaller dimensions, and most recently, Grinsztajn found a counterexample in dimension 63 \cite{Grinsztajn2026} working with GPT-5.5 Pro. OpenAI claims to have counterexamples for dimensions greater than 8 \cite{OpenAI2026Borsuk}.

For more details about the Borsuk partition problem we refer to the survey \cite{kalai2015some}.

\subsection{Reuleaux polyhedra}

Let $V\subset \mathbb{R}^3$ be a unit diameter finite set. The \emph{ball polyhedron} associated with $V$ is defined as $$\mathcal{B}(V)=\bigcap_{v\in V}B(v),$$
where $B(v)$ denotes the $3$-dimensional closed unit ball centered at $v$. Ball polyhedra have been widely studied; see \cite{bezdek2007ball} and \cite{kupitz2010ball} for further details.

The boundary points of a ball polyhedron can be naturally classified into faces, namely, \emph{spherical facets}, \emph{spherical edges}, and \emph{vertices}. The \emph{spherical facial structure} $\mathcal{SF}(\mathcal{B}(V))$ of a ball polyhedron $\mathcal{B}(V)$ is defined as the partial order of its faces under containment. A ball polyhedron $\mathcal{B}(V)$ is called \emph{standard} if $\mathcal{SF}(\mathcal{B}(V))$ is the face lattice of an abstract $3$-polytope (see \cite{abstract} for more information on abstract $3$-polytopes).

The set $V$ is called \emph{extremal} if the set of vertices of $\mathcal{B}(V)$, $\verti\mathcal{B}(V)$,  is exactly $V$. It is well known that a set $V$ is extremal if and only if its diameter graph has $2|V|-2$ edges ({Theorem 7.1, \cite{kupitz2010ball}). A standard ball polyhedron $\mathcal{B}(V)$ whose set $V$ is extremal is called a \emph{Reuleaux polyhedron}. Note that this definition differs from the one given originally by Sallee in \cite{sallee1970reuleaux} and from the one given by Hynd \cite{hynd2024density}, the latter being more flexible as it omits the standard condition.

We present recent results related to Reuleaux polyhedra and diameter graphs that are crucial for the proofs of Theorems \ref{Thm:Dolnikov} and \ref{Thm:Brass}. \\
\smallskip
\begin{comment}
    \begin{theorem}[Theorem 1 \cite{LOPEZCAMPOS2025104215}, 2025]\label{MainTheorem2025}
Let $V\subset \mathbb{R}^3$ be a finite set with unit diameter and $|V|=n \geq 4$. The following statements are equivalent
    \renewcommand{\theenumi}{\roman{enumi}}
    \begin{enumerate}
        \item $V$ has Borsuk number 4.
        \item There is a $V_1\subset V$ such that $\mathcal{B}(V_1)$ is a Reuleaux polyhedron.
    \end{enumerate}
\end{theorem}
\end{comment}

\begin{theorem}[Theorem 8 \cite{LOPEZCAMPOS2025104215}, 2025]\label{finiteC}
    Let $V\subset \mathbb{R}^3$ be a finite set of points with unit diameter and $|V| =n \geq 4$. The following three statements are equivalent:
    \begin{enumerate}[label=(\roman*)]
        \item\label{2eq} $\Bor (V)$ is 4-critical.
        \item\label{3eq} $\mathcal{B}(V)$ is a Reuleaux polyhedron.
    \end{enumerate}
\end{theorem}

\bigskip
 \begin{theorem}[Theorem 6  \cite{LOPEZCAMPOS2025104215}, 2025]\label{MeissnerThm2025}
     Let $\mathcal{B}(V)$ be a Reuleaux polyhedron. For every $\epsilon>0$ and every $v\in V$, there exist a closed covering $P_{\epsilon, v}=\{P_1, P_2, P_3, P_4\}$ of $\mathcal{B}(V)$, such that $v\in P_1$, $\diam (P_1)\leq \epsilon$ and $\diam (P_i) <1$, for $i=2,3,4$.
 \end{theorem}

\bigskip
 \begin{theorem} [Theorem 1.8 \cite{damasdi2026geometric}, 2026]\label{Thm:Damasdi}
     Let $G$ be a simple graph with at least one edge. Then the following are equivalent:
     \begin{enumerate}[label=\alph*)]
         \item $G$ is diameter graph in $\mathbb{R}^3$.
         \item $G$ is a subgraph of a simple non-bipartite quadrangulation of the projective plane.
     \end{enumerate}
 \end{theorem}

\section{Main contributions}\label{Sec:Main}

\subsection{On Dol'nikov's question}
In this subsection we prove Theorem \ref{Thm:Dolnikov} by proving its equivalent, Theorem \ref{Thm:DolnikovRephrase}. The outline of the proof is as follows:

\begin{enumerate}
    \item First, we define a family $\mathcal{F}$ of graphs (Definition \ref{Def: twisted graphs}) that are simple non-bipartite quadrangulations of the projective plane (Proposition \ref{prop:nbqp}).
    \item To make computations easier, we define a family $\mathcal{F}'$ whose members are subgraphs of some members of $\mathcal{F}$ (Definition \ref{Def:TwistedGrid}).
    \item Next, we compute the independence number of the elements of $\mathcal{F}'$ (Lemma \ref{lemma:alpha}).
    \item Finally, we show that for every $m\in\mathbb{N}$, there exist a graph $G(V,E)\in \mathcal{F}'$ such that its independence number $\alpha(G)<(|V|-m)/2$, and since $G$ has a geometric realization as a diameter graph by Theorem \ref{Thm:Damasdi}, the proof is complete.
\end{enumerate}

\bigskip
\begin{definition} \label{Def: twisted graphs}
For an odd integer $L\ge3$ and an even integer $W\geq 2$, we define the \emph{twisted prismoid} $Q_{L,\, W}$ as the graph with vertex set $V_{L,W}:=\{(x,y):0\le x\le L-1,\ 0\le y\le W\}\cup\{a\}$ (we refer to $a$ as the \emph{apex of} $Q_{L,\, W}$) and edges:
\begin{enumerate}
    \item \textbf{Vertical edges:} $(x,y)(x,y+1)$, for $0\leq x\leq L-1$ for $0\le y<W$ 
    \item \textbf{Horizontal edges:} $(x,y)(x+1,y)$, for $0\leq x< L-1$ for $0\leq y\leq W$.
    \item \textbf{Seam edges: } $(L-1,y)(0,W-y)$, for $0\leq y\leq W$.
    \item \textbf{Apex edges: } $a(x,0)$ if $x$ is even and $a(x,W)$ if $x$ is odd.
\end{enumerate}
Denote by $\mathcal{F}$ the family of twisted prismoid graphs.
\begin{figure}[h!]
    \centering
\begin{tikzpicture}[scale=1,>=stealth, every node/.style={font=\scriptsize}]
\draw[blue!70!black,thick] (0,0) -- (1,0);
\draw[blue!70!black,thick] (1,0) -- (2,0);
\draw[blue!70!black,thick] (2,0) -- (2.45,0);
\draw[blue!70!black,thick] (3.55,0) -- (4,0);
\node at (3.0,0) {$\cdots$};
\draw[blue!70!black,thick] (4,0) -- (5,0);
\draw[blue!70!black,thick] (0,1) -- (1,1);
\draw[blue!70!black,thick] (1,1) -- (2,1);
\draw[blue!70!black,thick] (2,1) -- (2.45,1);
\draw[blue!70!black,thick] (3.55,1) -- (4,1);
\node at (3.0,1) {$\cdots$};
\draw[blue!70!black,thick] (4,1) -- (5,1);
\draw[blue!70!black,thick] (0,3) -- (1,3);
\draw[blue!70!black,thick] (1,3) -- (2,3);
\draw[blue!70!black,thick] (2,3) -- (2.45,3);
\draw[blue!70!black,thick] (3.55,3) -- (4,3);
\node at (3.0,3) {$\cdots$};
\draw[blue!70!black,thick] (4,3) -- (5,3);
\draw[blue!70!black,thick] (0,5) -- (1,5);
\draw[blue!70!black,thick] (1,5) -- (2,5);
\draw[blue!70!black,thick] (2,5) -- (2.45,5);
\draw[blue!70!black,thick] (3.55,5) -- (4,5);
\node at (3.0,5) {$\cdots$};
\draw[blue!70!black,thick] (4,5) -- (5,5);
\draw[blue!70!black,thick] (0,6) -- (1,6);
\draw[blue!70!black,thick] (1,6) -- (2,6);
\draw[blue!70!black,thick] (2,6) -- (2.45,6);
\draw[blue!70!black,thick] (3.55,6) -- (4,6);
\node at (3.0,6) {$\cdots$};
\draw[blue!70!black,thick] (4,6) -- (5,6);
\draw[green!55!black,thick] (0,0) -- (0,1);
\draw[green!55!black,thick] (0,1) -- (0,1.45);
\draw[green!55!black,thick] (0,2.55) -- (0,3);
\node at (0,2.0) {$\vdots$};
\draw[green!55!black,thick] (0,3) -- (0,3.45);
\draw[green!55!black,thick] (0,4.55) -- (0,5);
\node at (0,4.0) {$\vdots$};
\draw[green!55!black,thick] (0,5) -- (0,6);
\draw[green!55!black,thick] (1,0) -- (1,1);
\draw[green!55!black,thick] (1,1) -- (1,1.45);
\draw[green!55!black,thick] (1,2.55) -- (1,3);
\node at (1,2.0) {$\vdots$};
\draw[green!55!black,thick] (1,3) -- (1,3.45);
\draw[green!55!black,thick] (1,4.55) -- (1,5);
\node at (1,4.0) {$\vdots$};
\draw[green!55!black,thick] (1,5) -- (1,6);
\draw[green!55!black,thick] (2,0) -- (2,1);
\draw[green!55!black,thick] (2,1) -- (2,1.45);
\draw[green!55!black,thick] (2,2.55) -- (2,3);
\node at (2,2.0) {$\vdots$};
\draw[green!55!black,thick] (2,3) -- (2,3.45);
\draw[green!55!black,thick] (2,4.55) -- (2,5);
\node at (2,4.0) {$\vdots$};
\draw[green!55!black,thick] (2,5) -- (2,6);
\draw[green!55!black,thick] (4,0) -- (4,1);
\draw[green!55!black,thick] (4,1) -- (4,1.45);
\draw[green!55!black,thick] (4,2.55) -- (4,3);
\node at (4,2.0) {$\vdots$};
\draw[green!55!black,thick] (4,3) -- (4,3.45);
\draw[green!55!black,thick] (4,4.55) -- (4,5);
\node at (4,4.0) {$\vdots$};
\draw[green!55!black,thick] (4,5) -- (4,6);
\draw[green!55!black,thick] (5,0) -- (5,1);
\draw[green!55!black,thick] (5,1) -- (5,1.45);
\draw[green!55!black,thick] (5,2.55) -- (5,3);
\node at (5,2.0) {$\vdots$};
\draw[green!55!black,thick] (5,3) -- (5,3.45);
\draw[green!55!black,thick] (5,4.55) -- (5,5);
\node at (5,4.0) {$\vdots$};
\draw[green!55!black,thick] (5,5) -- (5,6);
\node at (3,2) {$\ddots$};
\node at (3,4) {$\ddots$};
\draw[red!75!black,very thick] (5,0) -- (5.7,0);
\node[red!75!black,right] at (5.7,0) {$(0,W)$};
\draw[red!75!black,very thick] (0,0) -- (-0.7,0);
\node[red!75!black,left] at (-0.7,0) {$(L{-}1,W)$};
\draw[red!75!black,very thick] (5,1) -- (5.7,1);
\node[red!75!black,right] at (5.7,1) {$(0,W{-}1)$};
\draw[red!75!black,very thick] (0,1) -- (-0.7,1);
\node[red!75!black,left] at (-0.7,1) {$(L{-}1,W{-}1)$};
\draw[red!75!black,very thick] (5,3) -- (5.7,3);
\node[red!75!black,right] at (5.7,3.1) {$(0,\tfrac W2)$};
\draw[red!75!black,very thick] (0,3) -- (-0.7,3);
\node[red!75!black,left] at (-0.7,3) {$(L{-}1,\tfrac W2)$};
\draw[red!75!black,very thick] (5,5) -- (5.7,5);
\node[red!75!black,right] at (5.7,5) {$(0,1)$};
\draw[red!75!black,very thick] (0,5) -- (-0.7,5);
\node[red!75!black,left] at (-0.7,5) {$(L{-}1,1)$};
\draw[red!75!black,very thick] (5,6) -- (5.7,6);
\node[red!75!black,right] at (5.7,6) {$(0,0)$};
\draw[red!75!black,very thick] (0,6) -- (-0.7,6);
\node[red!75!black,left] at (-0.7,6) {$(L{-}1,0)$};
\draw[violet,thick,dashed] (0,0) -- (0,-0.6);
\node[violet,below] at (0,-0.6) {$a$};
\draw[violet,thick,dashed] (1,6) -- (1,6.6);
\node[violet,above] at (1,6.6) {$a$};
\draw[violet,thick,dashed] (2,0) -- (2,-0.6);
\node[violet,below] at (2,-0.6) {$a$};
\draw[violet,thick,dashed] (4,6) -- (4,6.6);
\node[violet,above] at (4,6.6) {$a$};
\draw[violet,thick,dashed] (5,0) -- (5,-0.6);
\node[violet,below] at (5,-0.6) {$a$};
\fill[black] (0,0) circle (2pt);
\fill[black] (0,1) circle (2pt);
\fill[orange!90!black] (0,3) circle (2.6pt);
\draw (0,3) circle (2.6pt);
\fill[black] (0,5) circle (2pt);
\fill[black] (0,6) circle (2pt);
\fill[black] (1,0) circle (2pt);
\fill[black] (1,1) circle (2pt);
\fill[orange!90!black] (1,3) circle (2.6pt);
\draw (1,3) circle (2.6pt);
\fill[black] (1,5) circle (2pt);
\fill[black] (1,6) circle (2pt);
\fill[black] (2,0) circle (2pt);
\fill[black] (2,1) circle (2pt);
\fill[orange!90!black] (2,3) circle (2.6pt);
\draw (2,3) circle (2.6pt);
\fill[black] (2,5) circle (2pt);
\fill[black] (2,6) circle (2pt);
\fill[black] (4,0) circle (2pt);
\fill[black] (4,1) circle (2pt);
\fill[orange!90!black] (4,3) circle (2.6pt);
\draw (4,3) circle (2.6pt);
\fill[black] (4,5) circle (2pt);
\fill[black] (4,6) circle (2pt);
\fill[black] (5,0) circle (2pt);
\fill[black] (5,1) circle (2pt);
\fill[orange!90!black] (5,3) circle (2.6pt);
\draw (5,3) circle (2.6pt);
\fill[black] (5,5) circle (2pt);
\fill[black] (5,6) circle (2pt);
\node[below left] at (0,0) {$(0,0)$};
\node[below left] at (1,0) {$(1,0)$};
\node[below left] at (2,0) {$(2,0)$};
\node[below left] at (4,0) {$(L{-}2,0)$};
\node[below right] at (5,0) {$(L{-}1,0)$};
\node[above left] at (0,6) {$(0,W)$};
\node[above left] at (1,6) {$(1,W)$};
\node[above left] at (2,6) {$(2,W)$};
\node[above left] at (4,6) {$(L{-}2,W)$};
\node[above right] at (5,6) {$(L{-}1,W)$};
\node[below left,orange!70!black] at (0,3) {$(0,\tfrac W2)$};
\node[below left,orange!70!black] at (1,3) {$(1,\tfrac W2)$};
\node[below left,orange!70!black] at (2,3) {$(2,\tfrac W2)$};
\node[below left,orange!70!black] at (4,3) {$(L{-}2,\tfrac W2)$};
\node[below right,orange!70!black] at (5,3) {$(L{-}1,\tfrac W2)$};
\node[orange!70!black,font=\small] at (-1.2,2.5) {middle row};
\end{tikzpicture}    \caption{Embedding of the twisted prismoid $Q_{L, W}$ in $\mathbb{P}^2$}
    \label{fig:QLW}
\end{figure}
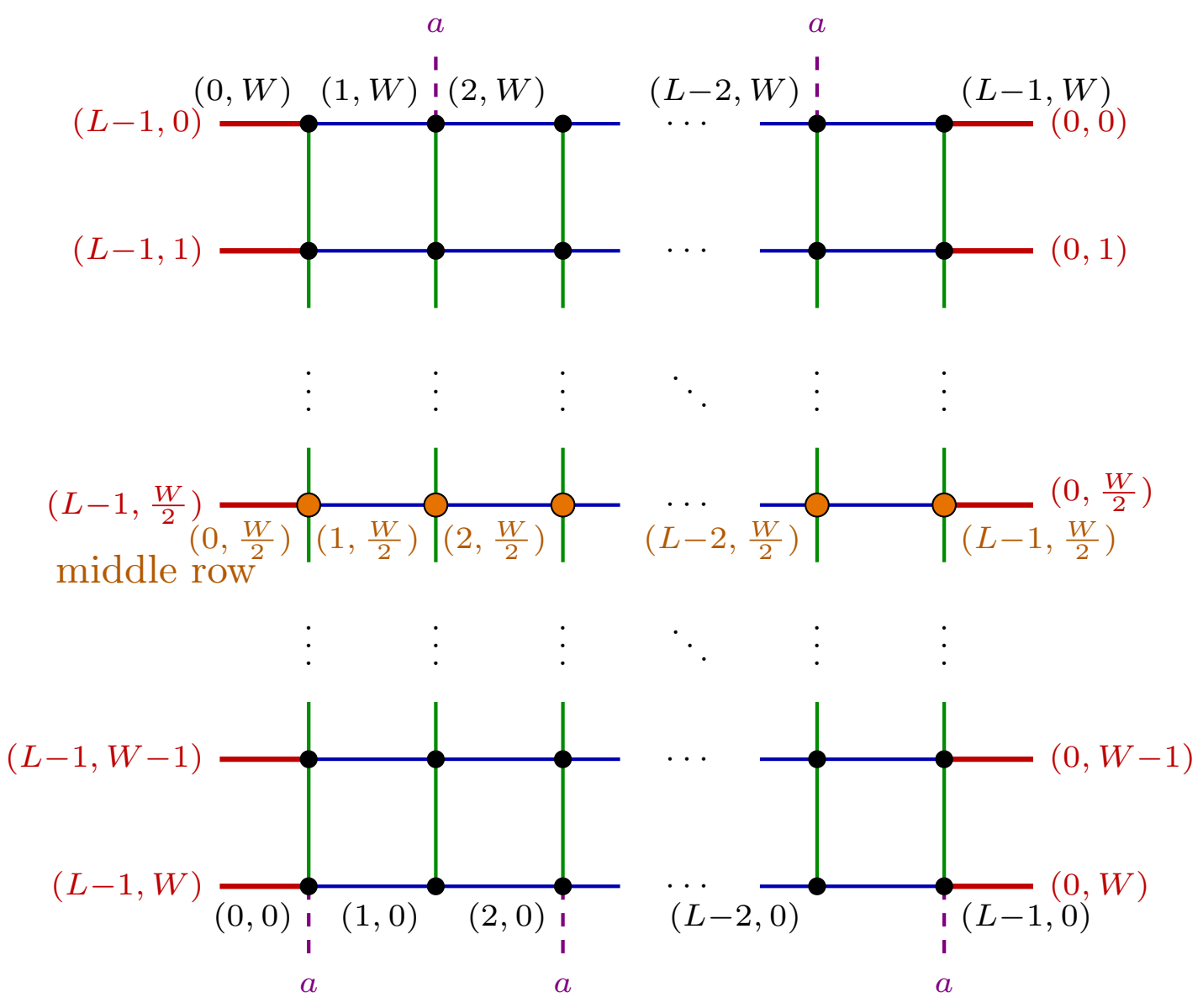
\end{definition}
\bigskip
\begin{proposition}\label{prop:nbqp}
$Q_{L,\, W}$ is a simple non-bipartite quadrangulation of the projective plane $\mathbb{P}^2$. That is, $Q_{L,\, W}$ has a $2$-cell embedding in $\mathbb{P}^2$ all of whose faces are bounded by $4$-cycles.
\end{proposition}
\begin{proof}
Take the rectangle $R=[0,L]\times[0,W]$ with the integer grid, and identify its left and right sides by $(0,t)\sim(L,W-t)$ (see Figure \ref{fig:QLW}). The result is a Möbius strip $M$, and the grid becomes an embedded graph on $M$: the vertices $(L,y)$ are identified with $(0,W-y)$, so the horizontal edges $(L-1,y)(L,y)$ become the seam edges $(L-1,y)(0,W-y)$, and the squares $[L-1,L]\times[y,y+1]$ become the seam faces $(L-1,y),(L-1,y+1),(0,W-y-1),(0,W-y)$. Therefore, all faces of this embedded graph are quadrangles so far. 

Note that the boundary of $M$ is the image of the top and bottom sides of $R$ on Figure.  Since the identification reverses orientation, the two sides form a single boundary cycle:
\[
b_0b_1\cdots b_{2L-1},\qquad b_j=(j,0),\ \ b_{L+j}=(j,W)\quad(0\le j<L),
\]
of length $2L$ (indeed $(L,0)\sim(0,W)$ and $(L,W)\sim(0,0)$). 

Next, we glue a disk $D$ along the boundary of $M$, turning it into the projective plane $\mathbb{P}^2=M\cup D$, and place the apex $a$ in the interior of $D$, joining it to $b_0,b_2,\dots,b_{2L-2}$ by disjoint arcs in $D$. Because $L$ is odd, the even-indexed boundary vertices are exactly the $(x,0)$ with $x$ even and the $(x,W)$ with $x$ odd, so these arcs are the apex edges of $Q_{L,\, W}$. 

Notice that these $L$ arcs divide $D$ into $L$ open disks, so the newly created faces are of the form $a\,b_{2i}b_{2i+1}b_{2i+2}$, for $0\leq i\leq L-1$ with indices $\mod 2L$, whose all are quadrangles. Thus every face of the embedding of $Q_{L,\, W}$ in $\mathbb{P}^2$ is a $4$-cycle, i.e.\ the embedding is a $2$-cell quadrangulation. 

Finally, the graph is simple (no loops or parallel edges arise from the identification because $L\ge3$ and $W\ge2$). It is non-bipartite because the middle row $(0,W/2),(1,W/2),\dots,(L-1,W/2)$ is a cycle of odd length $L$, whose last edge is the seam edge $(L-1,W/2)(0,W/2)$).
\end{proof}

\bigskip
\begin{remark}
      By Theorem \ref{Thm:Damasdi} and Proposition \ref{prop:nbqp}, there is a realization of $Q_{L,\, W}$ as a diameter graph in $\mathbb{R}^3$.  Furthermore, since $Q_{L,\, W}$ has $n:=L(W+1)+1$ vertices and $L(W+1)+LW+L=2L(W+1)=2n-2$ edges, this realization is extremal. In fact, it can be proved that these graphs are actually 4-critical, so they realize Reuleaux polyhedra by Theorem \ref{finiteC}. 
\end{remark}
\bigskip

\begin{definition}\label{Def:TwistedGrid}
    For an odd integer $L\geq3$, we define the \emph{twisted grid} $G_L$ as the induced subgraph of the twisted prismoid $Q_{L, 2L}$, with vertex set $V_L:=V_{L, 2L}\setminus\{a\}$.
    
Denote by $\mathcal{F}'$ the family of twisted grid graphs.
\end{definition}

Notice that $G_L$ can be embedded naturally on the Möbius grid as explained in the proof of Proposition \ref{prop:nbqp}, and its set of edges is the set of vertical, horizontal and seam edges of the twisted prismoid $Q_{L,2L}$. 
\bigskip

\begin{lemma}\label{lemma:alpha}
For an odd integer $L\geq3$, $\alpha(G_L)=L^2$.
\end{lemma}

\begin{proof}
For an easier computation of the independence number of $G_L$, we use the canonical double cover of the Möbius grid by the cylinder grid. Let $C_L$ be the cylinder grid with vertex set $V(C_L)=\mathbb{Z}_{2L}\times\{0,1,\dots,2L\}$ and edges $(\tilde x,y)(\tilde x{+}1,y)$, $(\tilde x,y)(\tilde x,y{+}1)$. Define the map $\pi:V(C_L)\to V_L$ as
\[
\pi(\tilde x,y)=\begin{cases}(\tilde x,y)&0\le\tilde x\le L-1,\\ (\tilde x-L,\ 2L-y)&L\le\tilde x\le 2L-1\end{cases}
\]
Notice that each vertex $(x,y)$ in $V_L$ has exactly two preimages in $V(C_L)$ under $\pi$, namely $(x,y)$ and $ (x+L, 2L-y)$. Furthermore, $\pi$ naturally sends 
the edges of $C_L$ onto the edges of $G_L$, having exactly two preimages each edge of $G_L$ as well. This is essentially the effect of the canonical double covering (see Figure \ref{Fig: Graph CL}).

\begin{figure}
    \centering
    \begin{tikzpicture}[scale=0.8, every node/.style={font=\scriptsize}]
\fill[blue!6] (-0.48,-1.0) rectangle (8.80,7.0); \fill[yellow!15] (8.80,-1.0) rectangle (18.24,7.0);
\node[font=\small,blue!60!black] at (4.00,7.3) {$\tilde x<L$:\ $\pi(\tilde x,y)=(\tilde x,y)$};
\node[font=\small,orange!70!black] at (13.60,7.3) {$\tilde x\ge L$:\ $\pi(\tilde x,y)=(\tilde x-L,\,2L-y)$};
\draw[blue!70!black,thick] (0.00,0) -- (1.60,0);
\draw[blue!70!black,thick] (1.60,0) -- (3.20,0);
\draw[blue!70!black,thick] (3.20,0) -- (3.80,0); \draw[blue!70!black,thick] (5.80,0) -- (6.40,0);
\node at (4.80,0) {$\cdots$};
\draw[blue!70!black,thick] (6.40,0) -- (8.00,0);
\draw[blue!70!black,thick] (0.00,1) -- (1.60,1);
\draw[blue!70!black,thick] (1.60,1) -- (3.20,1);
\draw[blue!70!black,thick] (3.20,1) -- (3.80,1); \draw[blue!70!black,thick] (5.80,1) -- (6.40,1);
\node at (4.80,1) {$\cdots$};
\draw[blue!70!black,thick] (6.40,1) -- (8.00,1);
\draw[blue!70!black,thick] (0.00,3) -- (1.60,3);
\draw[blue!70!black,thick] (1.60,3) -- (3.20,3);
\draw[blue!70!black,thick] (3.20,3) -- (3.80,3); \draw[blue!70!black,thick] (5.80,3) -- (6.40,3);
\node at (4.80,3) {$\cdots$};
\draw[blue!70!black,thick] (6.40,3) -- (8.00,3);
\draw[blue!70!black,thick] (0.00,5) -- (1.60,5);
\draw[blue!70!black,thick] (1.60,5) -- (3.20,5);
\draw[blue!70!black,thick] (3.20,5) -- (3.80,5); \draw[blue!70!black,thick] (5.80,5) -- (6.40,5);
\node at (4.80,5) {$\cdots$};
\draw[blue!70!black,thick] (6.40,5) -- (8.00,5);
\draw[blue!70!black,thick] (0.00,6) -- (1.60,6);
\draw[blue!70!black,thick] (1.60,6) -- (3.20,6);
\draw[blue!70!black,thick] (3.20,6) -- (3.80,6); \draw[blue!70!black,thick] (5.80,6) -- (6.40,6);
\node at (4.80,6) {$\cdots$};
\draw[blue!70!black,thick] (6.40,6) -- (8.00,6);
\draw[green!55!black,thick] (0.00,0) -- (0.00,1);
\draw[green!55!black,thick] (0.00,1) -- (0.00,1.45); \draw[green!55!black,thick] (0.00,2.55) -- (0.00,3);
\node at (0.00,2.0) {$\vdots$};
\draw[green!55!black,thick] (0.00,3) -- (0.00,3.45); \draw[green!55!black,thick] (0.00,4.55) -- (0.00,5);
\node at (0.00,4.0) {$\vdots$};
\draw[green!55!black,thick] (0.00,5) -- (0.00,6);
\draw[green!55!black,thick] (1.60,0) -- (1.60,1);
\draw[green!55!black,thick] (1.60,1) -- (1.60,1.45); \draw[green!55!black,thick] (1.60,2.55) -- (1.60,3);
\node at (1.60,2.0) {$\vdots$};
\draw[green!55!black,thick] (1.60,3) -- (1.60,3.45); \draw[green!55!black,thick] (1.60,4.55) -- (1.60,5);
\node at (1.60,4.0) {$\vdots$};
\draw[green!55!black,thick] (1.60,5) -- (1.60,6);
\draw[green!55!black,thick] (3.20,0) -- (3.20,1);
\draw[green!55!black,thick] (3.20,1) -- (3.20,1.45); \draw[green!55!black,thick] (3.20,2.55) -- (3.20,3);
\node at (3.20,2.0) {$\vdots$};
\draw[green!55!black,thick] (3.20,3) -- (3.20,3.45); \draw[green!55!black,thick] (3.20,4.55) -- (3.20,5);
\node at (3.20,4.0) {$\vdots$};
\draw[green!55!black,thick] (3.20,5) -- (3.20,6);
\draw[green!55!black,thick] (6.40,0) -- (6.40,1);
\draw[green!55!black,thick] (6.40,1) -- (6.40,1.45); \draw[green!55!black,thick] (6.40,2.55) -- (6.40,3);
\node at (6.40,2.0) {$\vdots$};
\draw[green!55!black,thick] (6.40,3) -- (6.40,3.45); \draw[green!55!black,thick] (6.40,4.55) -- (6.40,5);
\node at (6.40,4.0) {$\vdots$};
\draw[green!55!black,thick] (6.40,5) -- (6.40,6);
\draw[green!55!black,thick] (8.00,0) -- (8.00,1);
\draw[green!55!black,thick] (8.00,1) -- (8.00,1.45); \draw[green!55!black,thick] (8.00,2.55) -- (8.00,3);
\node at (8.00,2.0) {$\vdots$};
\draw[green!55!black,thick] (8.00,3) -- (8.00,3.45); \draw[green!55!black,thick] (8.00,4.55) -- (8.00,5);
\node at (8.00,4.0) {$\vdots$};
\draw[green!55!black,thick] (8.00,5) -- (8.00,6);
\node at (4.80,2) {$\ddots$};
\node at (4.80,4) {$\ddots$};
\draw[blue!70!black,thick] (9.60,0) -- (11.20,0);
\draw[blue!70!black,thick] (11.20,0) -- (12.80,0);
\draw[blue!70!black,thick] (12.80,0) -- (13.40,0); \draw[blue!70!black,thick] (15.40,0) -- (16.00,0);
\node at (14.40,0) {$\cdots$};
\draw[blue!70!black,thick] (16.00,0) -- (17.60,0);
\draw[blue!70!black,thick] (9.60,1) -- (11.20,1);
\draw[blue!70!black,thick] (11.20,1) -- (12.80,1);
\draw[blue!70!black,thick] (12.80,1) -- (13.40,1); \draw[blue!70!black,thick] (15.40,1) -- (16.00,1);
\node at (14.40,1) {$\cdots$};
\draw[blue!70!black,thick] (16.00,1) -- (17.60,1);
\draw[blue!70!black,thick] (9.60,3) -- (11.20,3);
\draw[blue!70!black,thick] (11.20,3) -- (12.80,3);
\draw[blue!70!black,thick] (12.80,3) -- (13.40,3); \draw[blue!70!black,thick] (15.40,3) -- (16.00,3);
\node at (14.40,3) {$\cdots$};
\draw[blue!70!black,thick] (16.00,3) -- (17.60,3);
\draw[blue!70!black,thick] (9.60,5) -- (11.20,5);
\draw[blue!70!black,thick] (11.20,5) -- (12.80,5);
\draw[blue!70!black,thick] (12.80,5) -- (13.40,5); \draw[blue!70!black,thick] (15.40,5) -- (16.00,5);
\node at (14.40,5) {$\cdots$};
\draw[blue!70!black,thick] (16.00,5) -- (17.60,5);
\draw[blue!70!black,thick] (9.60,6) -- (11.20,6);
\draw[blue!70!black,thick] (11.20,6) -- (12.80,6);
\draw[blue!70!black,thick] (12.80,6) -- (13.40,6); \draw[blue!70!black,thick] (15.40,6) -- (16.00,6);
\node at (14.40,6) {$\cdots$};
\draw[blue!70!black,thick] (16.00,6) -- (17.60,6);
\draw[green!55!black,thick] (9.60,0) -- (9.60,1);
\draw[green!55!black,thick] (9.60,1) -- (9.60,1.45); \draw[green!55!black,thick] (9.60,2.55) -- (9.60,3);
\node at (9.60,2.0) {$\vdots$};
\draw[green!55!black,thick] (9.60,3) -- (9.60,3.45); \draw[green!55!black,thick] (9.60,4.55) -- (9.60,5);
\node at (9.60,4.0) {$\vdots$};
\draw[green!55!black,thick] (9.60,5) -- (9.60,6);
\draw[green!55!black,thick] (11.20,0) -- (11.20,1);
\draw[green!55!black,thick] (11.20,1) -- (11.20,1.45); \draw[green!55!black,thick] (11.20,2.55) -- (11.20,3);
\node at (11.20,2.0) {$\vdots$};
\draw[green!55!black,thick] (11.20,3) -- (11.20,3.45); \draw[green!55!black,thick] (11.20,4.55) -- (11.20,5);
\node at (11.20,4.0) {$\vdots$};
\draw[green!55!black,thick] (11.20,5) -- (11.20,6);
\draw[green!55!black,thick] (12.80,0) -- (12.80,1);
\draw[green!55!black,thick] (12.80,1) -- (12.80,1.45); \draw[green!55!black,thick] (12.80,2.55) -- (12.80,3);
\node at (12.80,2.0) {$\vdots$};
\draw[green!55!black,thick] (12.80,3) -- (12.80,3.45); \draw[green!55!black,thick] (12.80,4.55) -- (12.80,5);
\node at (12.80,4.0) {$\vdots$};
\draw[green!55!black,thick] (12.80,5) -- (12.80,6);
\draw[green!55!black,thick] (16.00,0) -- (16.00,1);
\draw[green!55!black,thick] (16.00,1) -- (16.00,1.45); \draw[green!55!black,thick] (16.00,2.55) -- (16.00,3);
\node at (16.00,2.0) {$\vdots$};
\draw[green!55!black,thick] (16.00,3) -- (16.00,3.45); \draw[green!55!black,thick] (16.00,4.55) -- (16.00,5);
\node at (16.00,4.0) {$\vdots$};
\draw[green!55!black,thick] (16.00,5) -- (16.00,6);
\draw[green!55!black,thick] (17.60,0) -- (17.60,1);
\draw[green!55!black,thick] (17.60,1) -- (17.60,1.45); \draw[green!55!black,thick] (17.60,2.55) -- (17.60,3);
\node at (17.60,2.0) {$\vdots$};
\draw[green!55!black,thick] (17.60,3) -- (17.60,3.45); \draw[green!55!black,thick] (17.60,4.55) -- (17.60,5);
\node at (17.60,4.0) {$\vdots$};
\draw[green!55!black,thick] (17.60,5) -- (17.60,6);
\node at (14.40,2) {$\ddots$};
\node at (14.40,4) {$\ddots$};
\draw[blue!70!black,thick] (8.00,0) -- (9.60,0);
\draw[blue!70!black,thick] (17.60,0) -- (19.20,0);
\draw[fill=white] (19.20,0) circle (2pt); \node[right] at (19.20,0) {$=(0,0)$};
\draw[blue!70!black,thick] (8.00,1) -- (9.60,1);
\draw[blue!70!black,thick] (17.60,1) -- (19.20,1);
\draw[fill=white] (19.20,1) circle (2pt); \node[right] at (19.20,1) {$=(0,1)$};
\draw[blue!70!black,thick] (8.00,3) -- (9.60,3);
\draw[blue!70!black,thick] (17.60,3) -- (19.20,3);
\draw[fill=white] (19.20,3) circle (2pt); \node[right] at (19.20,3) {$=(0,L)$};
\draw[blue!70!black,thick] (8.00,5) -- (9.60,5);
\draw[blue!70!black,thick] (17.60,5) -- (19.20,5);
\draw[fill=white] (19.20,5) circle (2pt); \node[right] at (19.20,5) {$=(0,2L{-}1)$};
\draw[blue!70!black,thick] (8.00,6) -- (9.60,6);
\draw[blue!70!black,thick] (17.60,6) -- (19.20,6);
\draw[fill=white] (19.20,6) circle (2pt); \node[right] at (19.20,6) {$=(0,2L)$};
\draw[very thick,gray,->] (19.20,-0.4) -- (19.20,6.4); \draw[very thick,gray,->] (0.00,-0.4) -- (0.00,6.4);
\fill[black] (0.00,0) circle (2pt);
\fill[black] (0.00,1) circle (2pt);
\fill[orange!90!black] (0.00,3) circle (2.6pt);
\draw (0.00,3) circle (2.6pt);
\fill[black] (0.00,5) circle (2pt);
\fill[black] (0.00,6) circle (2pt);
\fill[black] (1.60,0) circle (2pt);
\fill[black] (1.60,1) circle (2pt);
\fill[orange!90!black] (1.60,3) circle (2.6pt);
\draw (1.60,3) circle (2.6pt);
\fill[black] (1.60,5) circle (2pt);
\fill[black] (1.60,6) circle (2pt);
\fill[black] (3.20,0) circle (2pt);
\fill[black] (3.20,1) circle (2pt);
\fill[orange!90!black] (3.20,3) circle (2.6pt);
\draw (3.20,3) circle (2.6pt);
\fill[black] (3.20,5) circle (2pt);
\fill[black] (3.20,6) circle (2pt);
\fill[black] (6.40,0) circle (2pt);
\fill[black] (6.40,1) circle (2pt);
\fill[orange!90!black] (6.40,3) circle (2.6pt);
\draw (6.40,3) circle (2.6pt);
\fill[black] (6.40,5) circle (2pt);
\fill[black] (6.40,6) circle (2pt);
\fill[black] (8.00,0) circle (2pt);
\fill[black] (8.00,1) circle (2pt);
\fill[orange!90!black] (8.00,3) circle (2.6pt);
\draw (8.00,3) circle (2.6pt);
\fill[black] (8.00,5) circle (2pt);
\fill[black] (8.00,6) circle (2pt);
\fill[black] (9.60,0) circle (2pt);
\fill[black] (9.60,1) circle (2pt);
\fill[orange!90!black] (9.60,3) circle (2.6pt);
\draw (9.60,3) circle (2.6pt);
\fill[black] (9.60,5) circle (2pt);
\fill[black] (9.60,6) circle (2pt);
\fill[black] (11.20,0) circle (2pt);
\fill[black] (11.20,1) circle (2pt);
\fill[orange!90!black] (11.20,3) circle (2.6pt);
\draw (11.20,3) circle (2.6pt);
\fill[black] (11.20,5) circle (2pt);
\fill[black] (11.20,6) circle (2pt);
\fill[black] (12.80,0) circle (2pt);
\fill[black] (12.80,1) circle (2pt);
\fill[orange!90!black] (12.80,3) circle (2.6pt);
\draw (12.80,3) circle (2.6pt);
\fill[black] (12.80,5) circle (2pt);
\fill[black] (12.80,6) circle (2pt);
\fill[black] (16.00,0) circle (2pt);
\fill[black] (16.00,1) circle (2pt);
\fill[orange!90!black] (16.00,3) circle (2.6pt);
\draw (16.00,3) circle (2.6pt);
\fill[black] (16.00,5) circle (2pt);
\fill[black] (16.00,6) circle (2pt);
\fill[black] (17.60,0) circle (2pt);
\fill[black] (17.60,1) circle (2pt);
\fill[orange!90!black] (17.60,3) circle (2.6pt);
\draw (17.60,3) circle (2.6pt);
\fill[black] (17.60,5) circle (2pt);
\fill[black] (17.60,6) circle (2pt);
\node[below] at (0.00,-0.05) {$(0,0)$};
\node[above] at (0.00,6.18) {$(0,2L)$};
\node[below,orange!70!black] at (0.00,2.95) {$(0,L)$};
\node[below] at (1.60,-0.05) {$(1,0)$};
\node[above] at (1.60,6.05) {$(1,2L)$};
\node[below,orange!70!black] at (1.60,2.95) {$(1,L)$};
\node[below] at (3.20,-0.05) {$(2,0)$};
\node[above] at (3.20,6.05) {$(2,2L)$};
\node[below,orange!70!black] at (3.20,2.95) {$(2,L)$};
\node[below] at (6.40,-0.05) {$(L{-}2,0)$};
\node[above] at (6.40,6.05) {$(L{-}2,2L)$};
\node[below,orange!70!black] at (6.40,2.95) {$(L{-}2,L)$};
\node[below] at (8.00,-0.05) {$(L{-}1,0)$};
\node[above] at (8.00,6.05) {$(L{-}1,2L)$};
\node[below,orange!70!black] at (8.00,2.95) {$(L{-}1,L)$};
\node[below] at (9.60,-0.05) {$(L,0)$};

\node[above] at (9.60,6.05) {$(L,2L)$};

\node[below,orange!70!black] at (9.60,2.95) {$(L,L)$};

\node[below] at (11.20,-0.05) {$(L{+}1,0)$};

\node[above] at (11.20,6.05) {$(L{+}1,2L)$};

\node[below,orange!70!black] at (11.20,2.95) {$(L{+}1,L)$};

\node[below] at (12.80,-0.05) {$(L{+}2,0)$};

\node[above] at (12.80,6.05) {$(L{+}2,2L)$};

\node[below,orange!70!black] at (12.80,2.95) {$(L{+}2,L)$};

\node[below] at (16.00,-0.05) {$(2L{-}2,0)$};

\node[above] at (16.00,6.05) {$(2L{-}2,2L)$};

\node[below,orange!70!black] at (16.00,2.95) {$(2L{-}2,L)$};

\node[below] at (17.60,-0.05) {$(2L{-}1,0)$};

\node[above] at (17.70,6.05) {$(2L{-}1,2L)$};

\node[below,orange!70!black] at (17.60,2.95) {$(2L{-}1,L)$};

\draw[purple,thick,<->,>=stealth] (0,0) to[bend left=22] node[pos=0.5,above,sloped,purple] {} (9.60,6.05);
\draw[purple,thick,<->,>=stealth] (8.10,1.12) to[bend left=22] node[pos=0.5,above,sloped,purple] {} (17.50,4.88);
\draw[purple,thick,<->,>=stealth] (3.20,3) to[bend left=22] node[pos=0.5,above,sloped,purple] {} (12.8,3);
\node[purple,font=\small] at (9.60,-1.4) {The red arrows denote the relation of the preimages of vertices of $G_L$  under $\pi$.};
\end{tikzpicture}
    \caption{The graph $C_L$ and the function $\pi$.}
    \label{Fig: Graph CL}
\end{figure}

\textbf{Claim 1:} Let $S\subset V_L$. Then $S$ is independent in $G_L$, if and only if $\St:=\pi^{-1}(S)$ is independent in $C_L$.

Indeed, since $\pi$ is a double covering, the subgraph of $C_L$ induced by $\St$ contains an edge if and only if the subgraph of $G_L$ induced by $S$ contains an edge.

By Claim 1, an upper bound on the size of an independent set $\St$ in $C_L$ yields an upper bound on the size of an independent set $S$ in $G_L$, since $|\St|=2|S|$. Then, we compute an upper bound for $|\St|$.

For $0\le y\le 2L$ let $\St_y=\St\cap(\Z_{2L}\times\{y\})$, let $e_y$ and $o_y$ be the numbers of its elements at even and at odd positions $\tilde x$ respectively, and define:
\[
d_y=L-e_y-o_y,\qquad \varphi_y=e_y-o_y .
\]
Notice that each row is a cycle of length $2L$, so $d_y\ge0$, and a maximum independent set of $C_L$ is a parity class, so
\begin{equation}\label{eq:full}
d_y=0\ \Longrightarrow\ |\varphi_y|=L .
\end{equation}
Summing over the $2L+1$ rows we obtain,
\begin{equation}\label{eq:count}
2|S|=|\St|=\sum_{y=0}^{2L}e_y+o_y=L(2L+1)-\sum_{y=0}^{2L}d_y=|V_L|-\sum_{y=0}^{2L}d_y .
\end{equation}

\textbf{Claim 2:} For $0\le y< 2L$
\begin{equation}\label{eq:consec}
|\varphi_y+\varphi_{y+1}|\le d_y+d_{y+1}
\end{equation}

Since vertices at the same column in consecutive rows are adjacent, there are at most $L$ elements in even-positions between $\St_y$ and $\St_{y+1}$, for $0\le y< 2L$. That is $e_y+e_{y+1}\le L$.  Likewise $o_y+o_{y+1}\le L$. 

On the other hand, by definition we have that $e_y+o_y=L-d_y$ and $e_{y+1} + o_{y+1} = L-d_{y+1}$.  Then
\begin{align*}
      L+(o_y + o_{y+1})&\geq (e_y + e_{y+1}) + (o_y + o_{y+1})\\ &= 2L - d_y - d_{y+1}\\
      \implies o_y + o_{y+1}&\geq L - d_y - d_{y+1},
\end{align*}
and 
\begin{align*}
      (e_y + e_{y+1})+L&\geq (e_y + e_{y+1}) + (o_y + o_{y+1})\\ &= 2L - d_y - d_{y+1}\\
      \implies (e_y + e_{y+1})&\geq L - d_y - d_{y+1}
\end{align*}
Therefore, 
\begin{align*}
    \varphi_y+\varphi_{y+1}&=(e_y+e_{y+1})-(o_y+o_{y+1})\\&\le L-(L-d_y-d_{y+1})\\&=d_y+d_{y+1},
\end{align*}
and
\begin{align*}
    \varphi_y+\varphi_{y+1}&=(e_y+e_{y+1})-(o_y+o_{y+1})\\&\ge (L-d_y-d_{y+1})-L\\&=-d_y-d_{y+1},
\end{align*}
This completes the proof of Claim 2.

Since $\pi(\tilde x,y)=\pi(\tilde x+L,2L-y)$ (with $\tilde x+L$ taken mod $2L$), we have $(\tilde x,y)\in\St$ if and only if $(\tilde x+L,2L-y)\in\St$. Because $L$ is odd and $2L$ is even, $\tilde x$ and $\tilde x+L$ have opposite parities, so $e_{2L-y}=o_y$ and $o_{2L-y}=e_y$. Hence, for each $0\le y\le 2L$,
\begin{align}
    \label{eq:sym1}\varphi_{2L-y}=e_{2L-y}-o_{2L-y}&=-o_y+e_y=-\varphi_y,\text{ and }\\\label{eq:sym2}
    d_{2L-y}=L-e_{2L-y}-o_{2L-y}&=L-o_y-e_y=d_y.
\end{align}

%This is where the twist enters: a checkerboard pattern would have $\varphi_y=\pm L$ with alternating signs, which is incompatible with \eqref{eq:sym} because $2L$ is even.

Let $\psi_j=(-1)^j\varphi_j$, for $j=0, 1,\dots,2L$. Then $|\psi_j-\psi_{j+1}|=|\varphi_j+\varphi_{j+1}|\le d_j+d_{j+1}$ by \eqref{eq:consec}, and $\psi_{2L-y}=(-1)^{2L-y}\varphi_{2L-y}=-(-1)^y\varphi_y=-\psi_y$ by \eqref{eq:sym1} and $2L$ even. Therefore, for $0\le y\le L$,
\[
2|\varphi_y|=|\psi_y-\psi_{2L-y}|\le\sum_{j=y}^{2L-y-1}|\psi_j-\psi_{j+1}|\le\sum_{j=y}^{2L-y-1}(d_j+d_{j+1})=2\sum_{j=y}^{2L-y}d_j-d_y-d_{2L-y},
\]
that is, by (\ref{eq:sym2}),
\begin{equation}\label{eq:tele}
\sum_{j=y}^{2L-y}d_j\ \ge\ |\varphi_y|+d_y.
\end{equation}

To bound $|S|$ from above, we proceed by cases. Recall that $d_y\geq 0$ for all $0\leq y\leq 2L$. If some row $y$ has $d_y=0$ (replacing $y$ by $2L-y$ we may assume $y\le L$), then $|\varphi_y|=L$ by \eqref{eq:full}, so \eqref{eq:tele} gives
\[
\sum_{i=0}^{2L}d_i\geq \sum_{j=y}^{2L-y}d_j\geq L.
\]
If no row has $d_y=0$, then 
\[
\sum_{i=0}^{2L}d_i\geq  2L+1>L.
\]
In both cases we have $\sum_i d_i\ge L$, and \eqref{eq:count} yields
\[
|S|=\frac{|V_L|-\sum_{i=0}^{2L}d_i }{2}\le\frac{|V_L|-L}{2}=\frac{L(2L+1)-L}{2}=L^2 .
\]

For the lower bound, we show an independent set in $G_L$ with $L^2$ vertices. Let
\[
A_1:=\{(x,y):\ y<L,\ x+y\ \text{odd}\}, \quad A_2:=\{(x,y):\ y>L,\ x+y\ \text{even}\},\quad  \text{and}\quad S_0=A_1\cup A_2.
\]
Observe that no horizontal or vertical edge of $G_L$ lies inside $A_1$ and $A_2$, and they are not vertically adjacent to each other because the middle row $y=L$ is empty. It remains to verify that $A_1$ and $A_2$ are not adjacent by a seam edge.

Consider a seam edge of the form $(L-1,y)(0,2L-y)$, where $0\leq y<L$. We have $(L-1,y)\in S_0$ if and only if $L-1+y$ is odd, which is equivalent to $y$ being odd. On the other hand, $(0,2L-y)\in S_0$ if and only if $2L-y$ is even (since $2L-y>L$), which is equivalent to $y$ being even. Thus, no seam edge of this form has both endpoints in $S_0$. The case $y>L$ follows by an analogous argument. Hence, $S_0$ is independent.

Finally, for each $0\leq y<L$, the rows $y$ and $2L-y$ contribute the vertices with $x+y$ odd and with $x+y$ even respectively. Together, these two rows contribute exactly $L$ vertices to $S_0$. Therefore, $|S_0|=L^2$.
\end{proof}

Now we proceed to prove Theorem \ref{Thm:DolnikovRephrase}.

\begin{proof}[Proof of Theorem \ref{Thm:DolnikovRephrase}]
    Let $m\in \mathbb{N}$ and choose an odd number $L>m$. Recall that the twisted grid $G_L$ is a subgraph of the twisted prismoid $Q_{L,2L}$, which, by Proposition \ref{prop:nbqp}, is a simple non-bipartite quadrangulation of the projective plane $\mathbb{P}^2$. Then, by Theorem \ref{Thm:Damasdi}, there exist a finite point set $S\subset \mathbb{R}^3$ such that $\Bor (S)\cong G_L$. Thus, by Lemma \ref{lemma:alpha}:
    \[
    \alpha(\Bor(S))= L^2=\frac{L(2L+1)-L}{2}<\frac{L(2L+1)-m}{2}=\frac{|S|-m}{2}\qedhere
    \]
\end{proof}
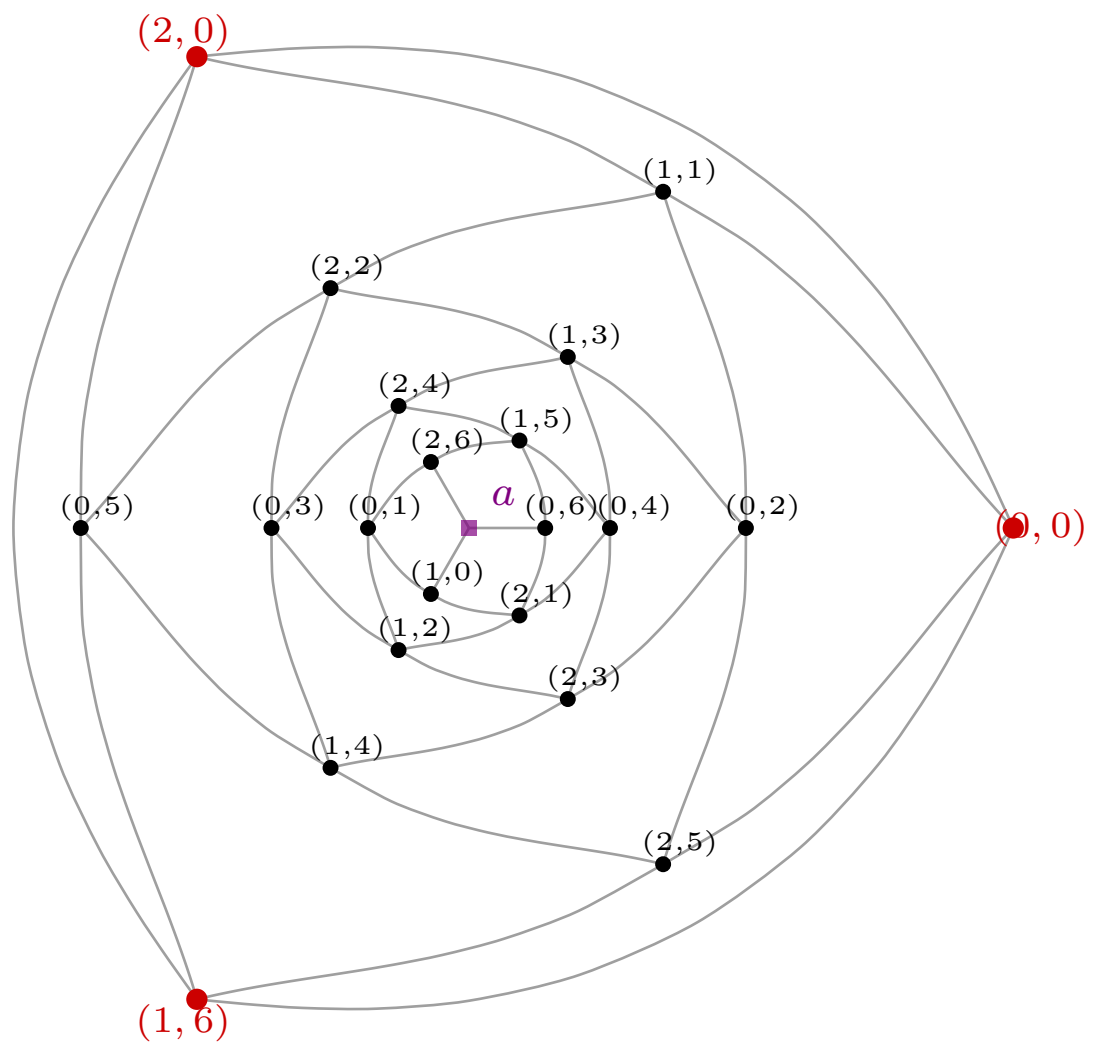
\begin{figure}[h!]
    \centering
    \begin{tikzpicture}[scale=0.5]
\tikzset{v/.style={circle,fill=black,inner sep=0pt,minimum size=3pt},o/.style={circle,fill=red!80!black,inner sep=0pt,minimum size=4pt},ap/.style={star,star points=5,fill=black,inner sep=0pt,minimum size=9pt}}
\draw[gray!75,line width=0.5pt] plot[smooth] coordinates {(7.200,0.000) (7.062,0.147) (6.913,0.312) (6.749,0.500) (6.569,0.714) (6.365,0.963) (6.130,1.255) (5.852,1.601) (5.509,2.017) (5.069,2.517) (4.480,3.113) (3.667,3.784) (2.568,4.448)};
\draw[gray!75,line width=0.5pt] plot[smooth] coordinates {(7.200,0.000) (6.940,-0.572) (6.598,-1.270) (6.119,-2.120) (5.416,-3.127) (4.388,-4.222) (3.014,-5.220) (1.463,-5.911) (-0.000,-6.254) (-1.224,-6.359) (-2.199,-6.349) (-2.974,-6.296) (-3.600,-6.235)};
\draw[gray!75,line width=0.5pt] plot[smooth] coordinates {(7.200,0.000) (6.940,0.572) (6.598,1.270) (6.119,2.120) (5.416,3.127) (4.388,4.222) (3.014,5.220) (1.463,5.911) (0.000,6.254) (-1.224,6.359) (-2.199,6.349) (-2.974,6.296) (-3.600,6.235)};
\draw[gray!75,line width=0.5pt] plot[smooth] coordinates {(7.200,0.000) (7.062,-0.147) (6.913,-0.312) (6.749,-0.500) (6.569,-0.714) (6.365,-0.963) (6.130,-1.255) (5.852,-1.601) (5.509,-2.017) (5.069,-2.517) (4.480,-3.113) (3.667,-3.784) (2.568,-4.448)};
\draw[gray!75,line width=0.5pt] plot[smooth] coordinates {(-1.336,0.000) (-1.311,-0.034) (-1.284,-0.072) (-1.255,-0.114) (-1.221,-0.163) (-1.183,-0.218) (-1.138,-0.282) (-1.084,-0.355) (-1.018,-0.441) (-0.934,-0.539) (-0.824,-0.649) (-0.682,-0.764) (-0.504,-0.873)};
\draw[gray!75,line width=0.5pt] plot[smooth] coordinates {(-1.336,0.000) (-1.330,-0.329) (-1.285,-0.597) (-1.227,-0.806) (-1.172,-0.971) (-1.122,-1.102) (-1.080,-1.210) (-1.044,-1.301) (-1.014,-1.379) (-0.988,-1.448) (-0.967,-1.509) (-0.948,-1.565) (-0.932,-1.615)};
\draw[gray!75,line width=0.5pt] plot[smooth] coordinates {(-1.336,0.000) (-1.330,0.329) (-1.285,0.597) (-1.227,0.806) (-1.172,0.971) (-1.122,1.102) (-1.080,1.210) (-1.044,1.301) (-1.014,1.379) (-0.988,1.448) (-0.967,1.509) (-0.948,1.565) (-0.932,1.615)};
\draw[gray!75,line width=0.5pt] plot[smooth] coordinates {(-1.336,0.000) (-1.311,0.034) (-1.284,0.072) (-1.255,0.114) (-1.221,0.163) (-1.183,0.218) (-1.138,0.282) (-1.084,0.355) (-1.018,0.441) (-0.934,0.539) (-0.824,0.649) (-0.682,0.764) (-0.504,0.873)};
\draw[gray!75,line width=0.5pt] plot[smooth] coordinates {(3.663,0.000) (3.646,0.916) (3.520,1.657) (3.363,2.234) (3.211,2.684) (3.076,3.044) (2.961,3.340) (2.865,3.589) (2.784,3.803) (2.716,3.991) (2.658,4.159) (2.610,4.310) (2.568,4.448)};
\draw[gray!75,line width=0.5pt] plot[smooth] coordinates {(3.663,0.000) (3.593,0.075) (3.517,0.159) (3.434,0.254) (3.342,0.363) (3.239,0.490) (3.119,0.638) (2.978,0.815) (2.803,1.026) (2.579,1.281) (2.279,1.584) (1.866,1.926) (1.307,2.263)};
\draw[gray!75,line width=0.5pt] plot[smooth] coordinates {(3.663,0.000) (3.593,-0.075) (3.517,-0.159) (3.434,-0.254) (3.342,-0.363) (3.239,-0.490) (3.119,-0.638) (2.978,-0.815) (2.803,-1.026) (2.579,-1.281) (2.279,-1.584) (1.866,-1.926) (1.307,-2.263)};
\draw[gray!75,line width=0.5pt] plot[smooth] coordinates {(3.663,0.000) (3.646,-0.916) (3.520,-1.657) (3.363,-2.234) (3.211,-2.684) (3.076,-3.044) (2.961,-3.340) (2.865,-3.589) (2.784,-3.803) (2.716,-3.991) (2.658,-4.159) (2.610,-4.310) (2.568,-4.448)};
\draw[gray!75,line width=0.5pt] plot[smooth] coordinates {(-2.613,0.000) (-2.563,-0.053) (-2.509,-0.114) (-2.450,-0.182) (-2.384,-0.260) (-2.310,-0.350) (-2.225,-0.456) (-2.124,-0.582) (-2.000,-0.733) (-1.840,-0.914) (-1.626,-1.131) (-1.331,-1.374) (-0.932,-1.615)};
\draw[gray!75,line width=0.5pt] plot[smooth] coordinates {(-2.613,0.000) (-2.601,-0.653) (-2.511,-1.182) (-2.399,-1.593) (-2.290,-1.915) (-2.194,-2.171) (-2.112,-2.382) (-2.044,-2.560) (-1.986,-2.713) (-1.937,-2.847) (-1.896,-2.966) (-1.861,-3.074) (-1.832,-3.173)};
\draw[gray!75,line width=0.5pt] plot[smooth] coordinates {(-2.613,0.000) (-2.601,0.653) (-2.511,1.182) (-2.399,1.593) (-2.290,1.915) (-2.194,2.171) (-2.112,2.382) (-2.044,2.560) (-1.986,2.713) (-1.937,2.847) (-1.896,2.966) (-1.861,3.074) (-1.832,3.173)};
\draw[gray!75,line width=0.5pt] plot[smooth] coordinates {(-2.613,0.000) (-2.563,0.053) (-2.509,0.114) (-2.450,0.182) (-2.384,0.260) (-2.310,0.350) (-2.225,0.456) (-2.124,0.582) (-2.000,0.733) (-1.840,0.914) (-1.626,1.131) (-1.331,1.374) (-0.932,1.615)};
\draw[gray!75,line width=0.5pt] plot[smooth] coordinates {(1.865,0.000) (1.856,0.466) (1.792,0.843) (1.712,1.136) (1.634,1.365) (1.566,1.549) (1.507,1.699) (1.458,1.826) (1.417,1.935) (1.382,2.031) (1.353,2.116) (1.328,2.193) (1.307,2.263)};
\draw[gray!75,line width=0.5pt] plot[smooth] coordinates {(1.865,0.000) (1.829,0.039) (1.790,0.082) (1.748,0.132) (1.701,0.188) (1.649,0.253) (1.588,0.330) (1.515,0.421) (1.426,0.529) (1.312,0.660) (1.159,0.814) (0.950,0.987) (0.668,1.157)};
\draw[gray!75,line width=0.5pt] plot[smooth] coordinates {(1.865,0.000) (1.829,-0.039) (1.790,-0.082) (1.748,-0.132) (1.701,-0.188) (1.649,-0.253) (1.588,-0.330) (1.515,-0.421) (1.426,-0.529) (1.312,-0.660) (1.159,-0.814) (0.950,-0.987) (0.668,-1.157)};
\draw[gray!75,line width=0.5pt] plot[smooth] coordinates {(1.865,0.000) (1.856,-0.466) (1.792,-0.843) (1.712,-1.136) (1.634,-1.365) (1.566,-1.549) (1.507,-1.699) (1.458,-1.826) (1.417,-1.935) (1.382,-2.031) (1.353,-2.116) (1.328,-2.193) (1.307,-2.263)};
\draw[gray!75,line width=0.5pt] plot[smooth] coordinates {(-5.136,0.000) (-5.037,-0.105) (-4.931,-0.223) (-4.814,-0.356) (-4.686,-0.510) (-4.540,-0.687) (-4.373,-0.895) (-4.174,-1.142) (-3.930,-1.438) (-3.616,-1.796) (-3.195,-2.220) (-2.616,-2.699) (-1.832,-3.173)};
\draw[gray!75,line width=0.5pt] plot[smooth] coordinates {(-5.136,0.000) (-5.111,-1.284) (-4.935,-2.323) (-4.715,-3.131) (-4.501,-3.763) (-4.312,-4.267) (-4.152,-4.682) (-4.016,-5.031) (-3.903,-5.332) (-3.808,-5.595) (-3.727,-5.830) (-3.658,-6.042) (-3.600,-6.235)};
\draw[gray!75,line width=0.5pt] plot[smooth] coordinates {(-5.136,0.000) (-5.111,1.284) (-4.935,2.323) (-4.715,3.131) (-4.501,3.763) (-4.312,4.267) (-4.152,4.682) (-4.016,5.031) (-3.903,5.332) (-3.808,5.595) (-3.727,5.830) (-3.658,6.042) (-3.600,6.235)};
\draw[gray!75,line width=0.5pt] plot[smooth] coordinates {(-5.136,0.000) (-5.037,0.105) (-4.931,0.223) (-4.814,0.356) (-4.686,0.510) (-4.540,0.687) (-4.373,0.895) (-4.174,1.142) (-3.930,1.438) (-3.616,1.796) (-3.195,2.220) (-2.616,2.699) (-1.832,3.173)};
\draw[gray!75,line width=0.5pt] plot[smooth] coordinates {(1.008,0.000) (1.003,0.208) (0.974,0.389) (0.934,0.539) (0.891,0.661) (0.850,0.761) (0.813,0.845) (0.780,0.915) (0.751,0.976) (0.726,1.029) (0.704,1.076) (0.685,1.119) (0.668,1.157)};
\draw[gray!75,line width=0.5pt] plot[smooth] coordinates {(1.008,0.000) (1.003,-0.208) (0.974,-0.389) (0.934,-0.539) (0.891,-0.661) (0.850,-0.761) (0.813,-0.845) (0.780,-0.915) (0.751,-0.976) (0.726,-1.029) (0.704,-1.076) (0.685,-1.119) (0.668,-1.157)};
\draw[gray!75,line width=0.5pt] plot[smooth] coordinates {(1.008,0.000) (0.986,0.000) (0.962,0.000) (0.936,0.000) (0.908,0.000) (0.878,0.000) (0.844,0.000) (0.805,0.000) (0.760,0.000) (0.706,0.000) (0.637,0.000) (0.533,0.000) (-0.000,0.000)};
\draw[gray!75,line width=0.5pt] plot[smooth] coordinates {(-0.504,-0.873) (-0.321,-0.973) (-0.150,-1.038) (-0.000,-1.078) (0.127,-1.102) (0.234,-1.117) (0.325,-1.127) (0.403,-1.133) (0.470,-1.139) (0.528,-1.144) (0.580,-1.148) (0.626,-1.153) (0.668,-1.157)};
\draw[gray!75,line width=0.5pt] plot[smooth] coordinates {(-0.504,-0.873) (-0.493,-0.854) (-0.481,-0.833) (-0.468,-0.811) (-0.454,-0.787) (-0.439,-0.760) (-0.422,-0.731) (-0.403,-0.697) (-0.380,-0.659) (-0.353,-0.612) (-0.318,-0.551) (-0.266,-0.461) (-0.000,0.000)};
\draw[gray!75,line width=0.5pt] plot[smooth] coordinates {(2.568,4.448) (1.444,5.068) (0.456,5.436) (-0.355,5.649) (-1.008,5.779) (-1.540,5.868) (-1.979,5.936) (-2.349,5.994) (-2.666,6.046) (-2.942,6.095) (-3.186,6.143) (-3.403,6.189) (-3.600,6.235)};
\draw[gray!75,line width=0.5pt] plot[smooth] coordinates {(2.568,4.448) (2.428,4.415) (2.272,4.382) (2.099,4.348) (1.902,4.313) (1.675,4.276) (1.412,4.235) (1.098,4.186) (0.719,4.122) (0.253,4.029) (-0.325,3.877) (-1.030,3.615) (-1.832,3.173)};
\draw[gray!75,line width=0.5pt] plot[smooth] coordinates {(-0.932,-1.615) (-0.881,-1.603) (-0.824,-1.592) (-0.760,-1.580) (-0.688,-1.568) (-0.605,-1.555) (-0.508,-1.540) (-0.393,-1.523) (-0.255,-1.500) (-0.085,-1.466) (0.125,-1.411) (0.380,-1.316) (0.668,-1.157)};
\draw[gray!75,line width=0.5pt] plot[smooth] coordinates {(-0.932,-1.615) (-0.525,-1.840) (-0.166,-1.973) (0.128,-2.051) (0.365,-2.098) (0.558,-2.130) (0.718,-2.155) (0.852,-2.176) (0.967,-2.195) (1.068,-2.212) (1.156,-2.230) (1.235,-2.246) (1.307,-2.263)};
\draw[gray!75,line width=0.5pt] plot[smooth] coordinates {(1.307,2.263) (0.735,2.579) (0.232,2.766) (-0.180,2.874) (-0.513,2.941) (-0.783,2.986) (-1.007,3.021) (-1.195,3.050) (-1.356,3.076) (-1.497,3.101) (-1.621,3.125) (-1.732,3.149) (-1.832,3.173)};
\draw[gray!75,line width=0.5pt] plot[smooth] coordinates {(1.307,2.263) (1.235,2.246) (1.156,2.230) (1.068,2.212) (0.967,2.195) (0.852,2.176) (0.718,2.155) (0.558,2.130) (0.365,2.098) (0.128,2.051) (-0.166,1.973) (-0.525,1.840) (-0.932,1.615)};
\draw[gray!75,line width=0.5pt] plot[smooth] coordinates {(-1.832,-3.173) (-1.732,-3.149) (-1.621,-3.125) (-1.497,-3.101) (-1.356,-3.076) (-1.195,-3.050) (-1.007,-3.021) (-0.783,-2.986) (-0.513,-2.941) (-0.180,-2.874) (0.232,-2.766) (0.735,-2.579) (1.307,-2.263)};
\draw[gray!75,line width=0.5pt] plot[smooth] coordinates {(-1.832,-3.173) (-1.030,-3.615) (-0.325,-3.877) (0.253,-4.029) (0.719,-4.122) (1.098,-4.186) (1.412,-4.235) (1.675,-4.276) (1.902,-4.313) (2.099,-4.348) (2.272,-4.382) (2.428,-4.415) (2.568,-4.448)};
\draw[gray!75,line width=0.5pt] plot[smooth] coordinates {(0.668,1.157) (0.380,1.316) (0.125,1.411) (-0.085,1.466) (-0.255,1.500) (-0.393,1.523) (-0.508,1.540) (-0.605,1.555) (-0.688,1.568) (-0.760,1.580) (-0.824,1.592) (-0.881,1.603) (-0.932,1.615)};
\draw[gray!75,line width=0.5pt] plot[smooth] coordinates {(0.668,1.157) (0.626,1.153) (0.580,1.148) (0.528,1.144) (0.470,1.139) (0.403,1.133) (0.325,1.127) (0.234,1.117) (0.127,1.102) (-0.000,1.078) (-0.150,1.038) (-0.321,0.973) (-0.504,0.873)};
\draw[gray!75,line width=0.5pt] plot[smooth] coordinates {(-3.600,-6.235) (-3.966,-5.724) (-4.399,-5.079) (-4.895,-4.239) (-5.416,-3.127) (-5.850,-1.689) (-6.027,0.000) (-5.850,1.689) (-5.416,3.127) (-4.895,4.239) (-4.399,5.079) (-3.966,5.724) (-3.600,6.235)};
\draw[gray!75,line width=0.5pt] plot[smooth] coordinates {(-3.600,-6.235) (-3.403,-6.189) (-3.186,-6.143) (-2.942,-6.095) (-2.666,-6.046) (-2.349,-5.994) (-1.979,-5.936) (-1.540,-5.868) (-1.008,-5.779) (-0.355,-5.649) (0.456,-5.436) (1.444,-5.068) (2.568,-4.448)};
\draw[gray!75,line width=0.5pt] plot[smooth] coordinates {(-0.504,0.873) (-0.493,0.854) (-0.481,0.833) (-0.468,0.811) (-0.454,0.787) (-0.439,0.760) (-0.422,0.731) (-0.403,0.697) (-0.380,0.659) (-0.353,0.612) (-0.318,0.551) (-0.266,0.461) (-0.000,0.000)};
\node[o] at (7.200,0.000) {}; \node[font=\scriptsize,red!80!black] at (7.560,0.000) {$(0,0)$};
\node[v] at (-1.336,0.000) {};
\node[font=\tiny,inner sep=0.5pt, text opacity=1] at (-1.116,0.270) {$\scriptscriptstyle(0,1)$};
\node[v] at (3.663,0.000) {};
\node[font=\tiny,inner sep=0.5pt, text opacity=1] at (3.883,0.270) {$\scriptscriptstyle(0,2)$};
\node[v] at (-2.613,0.000) {};
\node[font=\tiny,inner sep=0.5pt, text opacity=1] at (-2.393,0.270) {$\scriptscriptstyle(0,3)$};
\node[v] at (1.865,0.000) {};
\node[font=\tiny,inner sep=0.5pt, text opacity=1] at (2.185,0.270) {$\scriptscriptstyle(0,4)$};
\node[v] at (-5.136,0.000) {};
\node[font=\tiny,inner sep=0.5pt,text opacity=1] at (-4.916,0.270) {$\scriptscriptstyle(0,5)$};
\node[v] at (1.008,0.000) {};
\node[font=\tiny,inner sep=0.5pt, text opacity=1] at (1.228,0.270) {$\scriptscriptstyle(0,6)$};
\node[v] at (-0.504,-0.873) {};
\node[font=\tiny,inner sep=0.5pt, text opacity=1] at (-0.284,-0.603) {$\scriptscriptstyle(1,0)$};
\node[v] at (2.568,4.448) {};
\node[font=\tiny,inner sep=0.5pt, text opacity=1] at (2.788,4.718) {$\scriptscriptstyle(1,1)$};
\node[v] at (-0.932,-1.615) {};
\node[font=\tiny,inner sep=0.5pt, text opacity=1] at (-0.712,-1.345) {$\scriptscriptstyle(1,2)$};
\node[v] at (1.307,2.263) {};
\node[font=\tiny,inner sep=0.5pt, text opacity=1] at (1.527,2.533) {$\scriptscriptstyle(1,3)$};
\node[v] at (-1.832,-3.173) {};
\node[font=\tiny,inner sep=0.5pt, text opacity=1] at (-1.612,-2.903) {$\scriptscriptstyle(1,4)$};
\node[v] at (0.668,1.157) {};
\node[font=\tiny,inner sep=0.5pt, text opacity=1] at (0.888,1.427) {$\scriptscriptstyle(1,5)$};
\node[o] at (-3.600,-6.235) {}; \node[font=\scriptsize,red!80!black] at (-3.780,-6.547) {$(1,6)$};
\node[o] at (-3.600,6.235) {}; \node[font=\scriptsize,red!80!black] at (-3.780,6.547) {$(2,0)$};
\node[v] at (0.668,-1.157) {};
\node[font=\tiny,inner sep=0.5pt, text opacity=1] at (0.888,-0.887) {$\scriptscriptstyle(2,1)$};
\node[v] at (-1.832,3.173) {};
\node[font=\tiny,inner sep=0.5pt, text opacity=1] at (-1.612,3.443) {$\scriptscriptstyle(2,2)$};
\node[v] at (1.307,-2.263) {};
\node[font=\tiny,inner sep=0.5pt, text opacity=1] at (1.527,-1.993) {$\scriptscriptstyle(2,3)$};
\node[v] at (-0.932,1.615) {};
\node[font=\tiny,inner sep=0.5pt, text opacity=1] at (-0.712,1.885) {$\scriptscriptstyle(2,4)$};
\node[v] at (2.568,-4.448) {};
\node[font=\tiny,inner sep=0.5pt, text opacity=1] at (2.788,-4.178) {$\scriptscriptstyle(2,5)$};
\node[v] at (-0.504,0.873) {};
\node[font=\tiny,inner sep=0.5pt, text opacity=1] at (-0.284,1.143) {$\scriptscriptstyle(2,6)$};
\node[font=\tiny,inner sep=1.5pt,fill=violet,fill opacity=0.7,text opacity=1] at (-0.000,0.000) {}; \node[above right=1pt,font=\footnotesize, violet] at (-0.000,0.000) {$a$};
\end{tikzpicture}
    \caption{1-Skeleton of the Reuleaux polyhedron whose diameter graph is $Q_{3,\, 6}$}.
    \label{Fig:Q36}
\end{figure}

We would like to remark that the realization of $Q_{L,2L}$ (and consequently of $G_L$) as a diameter graph can be obtained without relying on Theorem~\ref{Thm:Damasdi}, by explicitly specifying the coordinates of its vertices for each $L$ (as Alon did for his construction in \cite{alon2003problems}). Figures \ref{Fig:Q36} and \ref{Fig:Q510} show the $1-$skeletons (the 0 and 1 dimensional faces) of the Reuleaux polyhedra whose diameter graphs are $Q_{3,\, 6}$ and $Q_{5,\, 10}$ respectively. Following this pattern, a symmetric embedding should not be difficult to find for every $Q_{L,2L}$, although doing so would require a significant amount of additional work.

\begin{figure}[h!]
    \centering
    \begin{tikzpicture}[scale=0.7]
\tikzset{v/.style={circle,fill=black,inner sep=0pt,minimum size=3pt},o/.style={circle,fill=red!80!black,inner sep=0pt,minimum size=4pt},ap/.style={star,star points=5,fill=black,inner sep=0pt,minimum size=9pt}}
\draw[gray!75,line width=0.5pt] plot[smooth] coordinates {(7.200,0.000) (7.092,0.186) (6.977,0.386) (6.852,0.601) (6.715,0.834) (6.565,1.085) (6.398,1.357) (6.210,1.651) (5.997,1.969) (5.754,2.311) (5.474,2.678) (5.149,3.065) (4.774,3.469)};
\draw[gray!75,line width=0.5pt] plot[smooth] coordinates {(7.200,0.000) (7.058,0.594) (6.868,1.230) (6.620,1.903) (6.307,2.598) (5.923,3.296) (5.472,3.976) (4.965,4.614) (4.419,5.195) (3.855,5.708) (3.292,6.152) (2.745,6.529) (2.225,6.848)};
\draw[gray!75,line width=0.5pt] plot[smooth] coordinates {(7.200,0.000) (7.058,-0.594) (6.868,-1.230) (6.620,-1.903) (6.307,-2.598) (5.923,-3.296) (5.472,-3.976) (4.965,-4.614) (4.419,-5.195) (3.855,-5.708) (3.292,-6.152) (2.745,-6.529) (2.225,-6.848)};
\draw[gray!75,line width=0.5pt] plot[smooth] coordinates {(7.200,0.000) (7.092,-0.186) (6.977,-0.386) (6.852,-0.601) (6.715,-0.834) (6.565,-1.085) (6.398,-1.357) (6.210,-1.651) (5.997,-1.969) (5.754,-2.311) (5.474,-2.678) (5.149,-3.065) (4.774,-3.469)};
\draw[gray!75,line width=0.5pt] plot[smooth] coordinates {(-1.209,0.000) (-1.192,-0.033) (-1.173,-0.068) (-1.153,-0.106) (-1.131,-0.147) (-1.106,-0.191) (-1.079,-0.238) (-1.048,-0.289) (-1.013,-0.343) (-0.973,-0.401) (-0.928,-0.462) (-0.875,-0.526) (-0.815,-0.592)};
\draw[gray!75,line width=0.5pt] plot[smooth] coordinates {(-1.209,0.000) (-1.222,-0.111) (-1.229,-0.213) (-1.231,-0.306) (-1.230,-0.392) (-1.226,-0.470) (-1.222,-0.541) (-1.216,-0.606) (-1.210,-0.666) (-1.204,-0.721) (-1.199,-0.772) (-1.193,-0.819) (-1.188,-0.863)};
\draw[gray!75,line width=0.5pt] plot[smooth] coordinates {(-1.209,0.000) (-1.222,0.111) (-1.229,0.213) (-1.231,0.306) (-1.230,0.392) (-1.226,0.470) (-1.222,0.541) (-1.216,0.606) (-1.210,0.666) (-1.204,0.721) (-1.199,0.772) (-1.193,0.819) (-1.188,0.863)};
\draw[gray!75,line width=0.5pt] plot[smooth] coordinates {(-1.209,0.000) (-1.192,0.033) (-1.173,0.068) (-1.153,0.106) (-1.131,0.147) (-1.106,0.191) (-1.079,0.238) (-1.048,0.289) (-1.013,0.343) (-0.973,0.401) (-0.928,0.462) (-0.875,0.526) (-0.815,0.592)};
\draw[gray!75,line width=0.5pt] plot[smooth] coordinates {(4.837,0.000) (4.891,0.448) (4.919,0.861) (4.929,1.239) (4.925,1.584) (4.913,1.897) (4.896,2.182) (4.876,2.443) (4.854,2.682) (4.833,2.902) (4.812,3.105) (4.793,3.293) (4.774,3.469)};
\draw[gray!75,line width=0.5pt] plot[smooth] coordinates {(4.837,0.000) (4.765,0.125) (4.687,0.259) (4.603,0.404) (4.511,0.560) (4.410,0.729) (4.298,0.911) (4.172,1.109) (4.029,1.323) (3.865,1.553) (3.677,1.799) (3.459,2.059) (3.207,2.330)};
\draw[gray!75,line width=0.5pt] plot[smooth] coordinates {(4.837,0.000) (4.765,-0.125) (4.687,-0.259) (4.603,-0.404) (4.511,-0.560) (4.410,-0.729) (4.298,-0.911) (4.172,-1.109) (4.029,-1.323) (3.865,-1.553) (3.677,-1.799) (3.459,-2.059) (3.207,-2.330)};
\draw[gray!75,line width=0.5pt] plot[smooth] coordinates {(4.837,0.000) (4.891,-0.448) (4.919,-0.861) (4.929,-1.239) (4.925,-1.584) (4.913,-1.897) (4.896,-2.182) (4.876,-2.443) (4.854,-2.682) (4.833,-2.902) (4.812,-3.105) (4.793,-3.293) (4.774,-3.469)};
\draw[gray!75,line width=0.5pt] plot[smooth] coordinates {(-1.790,0.000) (-1.763,-0.046) (-1.735,-0.096) (-1.704,-0.150) (-1.670,-0.208) (-1.632,-0.271) (-1.591,-0.338) (-1.544,-0.412) (-1.492,-0.491) (-1.431,-0.576) (-1.362,-0.667) (-1.281,-0.763) (-1.188,-0.863)};
\draw[gray!75,line width=0.5pt] plot[smooth] coordinates {(-1.790,0.000) (-1.810,-0.166) (-1.820,-0.319) (-1.824,-0.458) (-1.822,-0.586) (-1.818,-0.702) (-1.812,-0.807) (-1.804,-0.904) (-1.796,-0.992) (-1.788,-1.074) (-1.780,-1.149) (-1.773,-1.218) (-1.766,-1.283)};
\draw[gray!75,line width=0.5pt] plot[smooth] coordinates {(-1.790,0.000) (-1.810,0.166) (-1.820,0.319) (-1.824,0.458) (-1.822,0.586) (-1.818,0.702) (-1.812,0.807) (-1.804,0.904) (-1.796,0.992) (-1.788,1.074) (-1.780,1.149) (-1.773,1.218) (-1.766,1.283)};
\draw[gray!75,line width=0.5pt] plot[smooth] coordinates {(-1.790,0.000) (-1.763,0.046) (-1.735,0.096) (-1.704,0.150) (-1.670,0.208) (-1.632,0.271) (-1.591,0.338) (-1.544,0.412) (-1.492,0.491) (-1.431,0.576) (-1.362,0.667) (-1.281,0.763) (-1.188,0.863)};
\draw[gray!75,line width=0.5pt] plot[smooth] coordinates {(3.249,0.000) (3.286,0.301) (3.305,0.579) (3.311,0.833) (3.309,1.064) (3.301,1.274) (3.289,1.466) (3.276,1.641) (3.261,1.802) (3.247,1.950) (3.233,2.086) (3.220,2.213) (3.207,2.330)};
\draw[gray!75,line width=0.5pt] plot[smooth] coordinates {(3.249,0.000) (3.201,0.084) (3.149,0.174) (3.092,0.271) (3.031,0.376) (2.963,0.490) (2.887,0.612) (2.803,0.745) (2.707,0.889) (2.597,1.043) (2.470,1.208) (2.324,1.383) (2.155,1.566)};
\draw[gray!75,line width=0.5pt] plot[smooth] coordinates {(3.249,0.000) (3.201,-0.084) (3.149,-0.174) (3.092,-0.271) (3.031,-0.376) (2.963,-0.490) (2.887,-0.612) (2.803,-0.745) (2.707,-0.889) (2.597,-1.043) (2.470,-1.208) (2.324,-1.383) (2.155,-1.566)};
\draw[gray!75,line width=0.5pt] plot[smooth] coordinates {(3.249,0.000) (3.286,-0.301) (3.305,-0.579) (3.311,-0.833) (3.309,-1.064) (3.301,-1.274) (3.289,-1.466) (3.276,-1.641) (3.261,-1.802) (3.247,-1.950) (3.233,-2.086) (3.220,-2.213) (3.207,-2.330)};
\draw[gray!75,line width=0.5pt] plot[smooth] coordinates {(-2.663,0.000) (-2.624,-0.069) (-2.581,-0.143) (-2.535,-0.223) (-2.484,-0.308) (-2.428,-0.401) (-2.367,-0.502) (-2.297,-0.611) (-2.219,-0.728) (-2.129,-0.855) (-2.025,-0.991) (-1.905,-1.134) (-1.766,-1.283)};
\draw[gray!75,line width=0.5pt] plot[smooth] coordinates {(-2.663,0.000) (-2.693,-0.247) (-2.709,-0.474) (-2.714,-0.682) (-2.712,-0.872) (-2.705,-1.045) (-2.696,-1.202) (-2.685,-1.345) (-2.673,-1.477) (-2.661,-1.598) (-2.650,-1.710) (-2.639,-1.813) (-2.629,-1.910)};
\draw[gray!75,line width=0.5pt] plot[smooth] coordinates {(-2.663,0.000) (-2.693,0.247) (-2.709,0.474) (-2.714,0.682) (-2.712,0.872) (-2.705,1.045) (-2.696,1.202) (-2.685,1.345) (-2.673,1.477) (-2.661,1.598) (-2.650,1.710) (-2.639,1.813) (-2.629,1.910)};
\draw[gray!75,line width=0.5pt] plot[smooth] coordinates {(-2.663,0.000) (-2.624,0.069) (-2.581,0.143) (-2.535,0.223) (-2.484,0.308) (-2.428,0.401) (-2.367,0.502) (-2.297,0.611) (-2.219,0.728) (-2.129,0.855) (-2.025,0.991) (-1.905,1.134) (-1.766,1.283)};
\draw[gray!75,line width=0.5pt] plot[smooth] coordinates {(2.183,0.000) (2.208,0.202) (2.220,0.389) (2.225,0.559) (2.223,0.715) (2.218,0.856) (2.210,0.985) (2.201,1.103) (2.191,1.211) (2.181,1.310) (2.172,1.401) (2.163,1.486) (2.155,1.566)};
\draw[gray!75,line width=0.5pt] plot[smooth] coordinates {(2.183,0.000) (2.151,0.056) (2.116,0.117) (2.078,0.183) (2.036,0.253) (1.991,0.329) (1.940,0.412) (1.883,0.501) (1.819,0.597) (1.745,0.701) (1.660,0.812) (1.562,0.930) (1.448,1.052)};
\draw[gray!75,line width=0.5pt] plot[smooth] coordinates {(2.183,0.000) (2.151,-0.056) (2.116,-0.117) (2.078,-0.183) (2.036,-0.253) (1.991,-0.329) (1.940,-0.412) (1.883,-0.501) (1.819,-0.597) (1.745,-0.701) (1.660,-0.812) (1.562,-0.930) (1.448,-1.052)};
\draw[gray!75,line width=0.5pt] plot[smooth] coordinates {(2.183,0.000) (2.208,-0.202) (2.220,-0.389) (2.225,-0.559) (2.223,-0.715) (2.218,-0.856) (2.210,-0.985) (2.201,-1.103) (2.191,-1.211) (2.181,-1.310) (2.172,-1.401) (2.163,-1.486) (2.155,-1.566)};
\draw[gray!75,line width=0.5pt] plot[smooth] coordinates {(-3.965,0.000) (-3.905,-0.102) (-3.842,-0.213) (-3.773,-0.331) (-3.698,-0.459) (-3.615,-0.597) (-3.523,-0.747) (-3.419,-0.909) (-3.302,-1.084) (-3.168,-1.273) (-3.014,-1.474) (-2.835,-1.688) (-2.629,-1.910)};
\draw[gray!75,line width=0.5pt] plot[smooth] coordinates {(-3.965,0.000) (-4.009,-0.367) (-4.032,-0.706) (-4.040,-1.016) (-4.037,-1.298) (-4.027,-1.555) (-4.013,-1.789) (-3.996,-2.003) (-3.979,-2.198) (-3.961,-2.379) (-3.944,-2.545) (-3.928,-2.699) (-3.913,-2.843)};
\draw[gray!75,line width=0.5pt] plot[smooth] coordinates {(-3.965,0.000) (-4.009,0.367) (-4.032,0.706) (-4.040,1.016) (-4.037,1.298) (-4.027,1.555) (-4.013,1.789) (-3.996,2.003) (-3.979,2.198) (-3.961,2.379) (-3.944,2.545) (-3.928,2.699) (-3.913,2.843)};
\draw[gray!75,line width=0.5pt] plot[smooth] coordinates {(-3.965,0.000) (-3.905,0.102) (-3.842,0.213) (-3.773,0.331) (-3.698,0.459) (-3.615,0.597) (-3.523,0.747) (-3.419,0.909) (-3.302,1.084) (-3.168,1.273) (-3.014,1.474) (-2.835,1.688) (-2.629,1.910)};
\draw[gray!75,line width=0.5pt] plot[smooth] coordinates {(1.469,0.000) (1.485,0.136) (1.494,0.261) (1.496,0.375) (1.495,0.480) (1.491,0.575) (1.486,0.661) (1.480,0.741) (1.473,0.813) (1.466,0.880) (1.460,0.942) (1.454,0.999) (1.448,1.052)};
\draw[gray!75,line width=0.5pt] plot[smooth] coordinates {(1.469,0.000) (1.447,0.039) (1.424,0.080) (1.398,0.124) (1.371,0.172) (1.340,0.224) (1.306,0.280) (1.268,0.340) (1.225,0.406) (1.176,0.476) (1.119,0.550) (1.054,0.629) (0.978,0.711)};
\draw[gray!75,line width=0.5pt] plot[smooth] coordinates {(1.469,0.000) (1.447,-0.039) (1.424,-0.080) (1.398,-0.124) (1.371,-0.172) (1.340,-0.224) (1.306,-0.280) (1.268,-0.340) (1.225,-0.406) (1.176,-0.476) (1.119,-0.550) (1.054,-0.629) (0.978,-0.711)};
\draw[gray!75,line width=0.5pt] plot[smooth] coordinates {(1.469,0.000) (1.485,-0.136) (1.494,-0.261) (1.496,-0.375) (1.495,-0.480) (1.491,-0.575) (1.486,-0.661) (1.480,-0.741) (1.473,-0.813) (1.466,-0.880) (1.460,-0.942) (1.454,-0.999) (1.448,-1.052)};
\draw[gray!75,line width=0.5pt] plot[smooth] coordinates {(-5.901,0.000) (-5.813,-0.153) (-5.718,-0.316) (-5.616,-0.493) (-5.504,-0.683) (-5.381,-0.889) (-5.244,-1.112) (-5.090,-1.353) (-4.915,-1.614) (-4.716,-1.894) (-4.486,-2.195) (-4.221,-2.512) (-3.913,-2.843)};
\draw[gray!75,line width=0.5pt] plot[smooth] coordinates {(-5.901,0.000) (-5.968,-0.547) (-6.002,-1.051) (-6.014,-1.512) (-6.009,-1.932) (-5.994,-2.314) (-5.973,-2.663) (-5.949,-2.981) (-5.923,-3.273) (-5.897,-3.541) (-5.871,-3.789) (-5.847,-4.018) (-5.825,-4.232)};
\draw[gray!75,line width=0.5pt] plot[smooth] coordinates {(-5.901,0.000) (-5.968,0.547) (-6.002,1.051) (-6.014,1.512) (-6.009,1.932) (-5.994,2.314) (-5.973,2.663) (-5.949,2.981) (-5.923,3.273) (-5.897,3.541) (-5.871,3.789) (-5.847,4.018) (-5.825,4.232)};
\draw[gray!75,line width=0.5pt] plot[smooth] coordinates {(-5.901,0.000) (-5.813,0.153) (-5.718,0.316) (-5.616,0.493) (-5.504,0.683) (-5.381,0.889) (-5.244,1.112) (-5.090,1.353) (-4.915,1.614) (-4.716,1.894) (-4.486,2.195) (-4.221,2.512) (-3.913,2.843)};
\draw[gray!75,line width=0.5pt] plot[smooth] coordinates {(1.008,0.000) (1.017,0.089) (1.022,0.171) (1.023,0.248) (1.021,0.318) (1.017,0.382) (1.013,0.441) (1.007,0.496) (1.001,0.546) (0.995,0.592) (0.989,0.634) (0.984,0.674) (0.978,0.711)};
\draw[gray!75,line width=0.5pt] plot[smooth] coordinates {(1.008,0.000) (1.017,-0.089) (1.022,-0.171) (1.023,-0.248) (1.021,-0.318) (1.017,-0.382) (1.013,-0.441) (1.007,-0.496) (1.001,-0.546) (0.995,-0.592) (0.989,-0.634) (0.984,-0.674) (0.978,-0.711)};
\draw[gray!75,line width=0.5pt] plot[smooth] coordinates {(1.008,0.000) (0.982,0.000) (0.955,0.000) (0.926,0.000) (0.894,0.000) (0.860,0.000) (0.822,0.000) (0.779,0.000) (0.729,0.000) (0.670,0.000) (0.594,0.000) (0.484,0.000) (-0.000,0.000)};
\draw[gray!75,line width=0.5pt] plot[smooth] coordinates {(-0.815,-0.592) (-0.771,-0.670) (-0.726,-0.739) (-0.682,-0.802) (-0.639,-0.857) (-0.598,-0.907) (-0.560,-0.952) (-0.523,-0.993) (-0.489,-1.030) (-0.457,-1.064) (-0.428,-1.095) (-0.400,-1.123) (-0.374,-1.150)};
\draw[gray!75,line width=0.5pt] plot[smooth] coordinates {(-0.815,-0.592) (-0.795,-0.577) (-0.773,-0.561) (-0.749,-0.544) (-0.724,-0.526) (-0.696,-0.505) (-0.665,-0.483) (-0.630,-0.458) (-0.590,-0.428) (-0.542,-0.394) (-0.481,-0.349) (-0.392,-0.285) (-0.000,0.000)};
\draw[gray!75,line width=0.5pt] plot[smooth] coordinates {(4.774,3.469) (4.506,3.950) (4.238,4.378) (3.976,4.758) (3.726,5.095) (3.489,5.396) (3.267,5.665) (3.060,5.908) (2.868,6.129) (2.689,6.330) (2.523,6.516) (2.369,6.688) (2.225,6.848)};
\draw[gray!75,line width=0.5pt] plot[smooth] coordinates {(4.774,3.469) (4.613,3.540) (4.440,3.617) (4.254,3.700) (4.051,3.788) (3.830,3.882) (3.589,3.982) (3.322,4.086) (3.028,4.195) (2.702,4.305) (2.339,4.412) (1.938,4.513) (1.495,4.600)};
\draw[gray!75,line width=0.5pt] plot[smooth] coordinates {(-1.188,-0.863) (-1.148,-0.882) (-1.105,-0.902) (-1.058,-0.923) (-1.008,-0.945) (-0.953,-0.969) (-0.892,-0.994) (-0.826,-1.021) (-0.753,-1.048) (-0.672,-1.076) (-0.582,-1.103) (-0.483,-1.128) (-0.374,-1.150)};
\draw[gray!75,line width=0.5pt] plot[smooth] coordinates {(-1.188,-0.863) (-1.122,-0.983) (-1.055,-1.089) (-0.990,-1.183) (-0.928,-1.267) (-0.869,-1.342) (-0.813,-1.408) (-0.762,-1.469) (-0.714,-1.524) (-0.669,-1.574) (-0.628,-1.620) (-0.589,-1.663) (-0.553,-1.702)};
\draw[gray!75,line width=0.5pt] plot[smooth] coordinates {(3.207,2.330) (3.027,2.654) (2.847,2.941) (2.671,3.196) (2.503,3.423) (2.344,3.625) (2.195,3.806) (2.056,3.969) (1.927,4.117) (1.807,4.253) (1.695,4.377) (1.591,4.493) (1.495,4.600)};
\draw[gray!75,line width=0.5pt] plot[smooth] coordinates {(3.207,2.330) (3.099,2.378) (2.983,2.430) (2.858,2.485) (2.722,2.545) (2.573,2.608) (2.411,2.675) (2.232,2.745) (2.034,2.818) (1.815,2.892) (1.572,2.964) (1.302,3.032) (1.004,3.090)};
\draw[gray!75,line width=0.5pt] plot[smooth] coordinates {(-1.766,-1.283) (-1.707,-1.310) (-1.643,-1.338) (-1.574,-1.369) (-1.499,-1.402) (-1.417,-1.436) (-1.327,-1.473) (-1.229,-1.512) (-1.120,-1.552) (-0.999,-1.593) (-0.865,-1.633) (-0.717,-1.670) (-0.553,-1.702)};
\draw[gray!75,line width=0.5pt] plot[smooth] coordinates {(-1.766,-1.283) (-1.667,-1.461) (-1.568,-1.620) (-1.471,-1.760) (-1.378,-1.885) (-1.291,-1.996) (-1.209,-2.096) (-1.132,-2.186) (-1.061,-2.267) (-0.995,-2.342) (-0.933,-2.410) (-0.876,-2.474) (-0.823,-2.533)};
\draw[gray!75,line width=0.5pt] plot[smooth] coordinates {(2.155,1.566) (2.034,1.783) (1.913,1.976) (1.795,2.147) (1.682,2.300) (1.575,2.435) (1.475,2.557) (1.381,2.666) (1.294,2.766) (1.214,2.857) (1.139,2.941) (1.069,3.018) (1.004,3.090)};
\draw[gray!75,line width=0.5pt] plot[smooth] coordinates {(2.155,1.566) (2.082,1.598) (2.004,1.633) (1.920,1.670) (1.828,1.710) (1.729,1.752) (1.620,1.797) (1.499,1.844) (1.367,1.893) (1.219,1.943) (1.056,1.992) (0.875,2.037) (0.675,2.076)};
\draw[gray!75,line width=0.5pt] plot[smooth] coordinates {(-2.629,-1.910) (-2.540,-1.949) (-2.445,-1.992) (-2.342,-2.037) (-2.231,-2.086) (-2.109,-2.138) (-1.976,-2.193) (-1.829,-2.250) (-1.667,-2.310) (-1.488,-2.370) (-1.288,-2.430) (-1.067,-2.485) (-0.823,-2.533)};
\draw[gray!75,line width=0.5pt] plot[smooth] coordinates {(-2.629,-1.910) (-2.481,-2.175) (-2.334,-2.411) (-2.189,-2.620) (-2.052,-2.806) (-1.921,-2.971) (-1.799,-3.119) (-1.685,-3.253) (-1.579,-3.375) (-1.481,-3.486) (-1.389,-3.588) (-1.304,-3.682) (-1.225,-3.770)};
\draw[gray!75,line width=0.5pt] plot[smooth] coordinates {(1.448,1.052) (1.367,1.198) (1.286,1.328) (1.206,1.443) (1.130,1.545) (1.058,1.636) (0.991,1.718) (0.928,1.792) (0.870,1.858) (0.816,1.920) (0.765,1.976) (0.718,2.028) (0.675,2.076)};
\draw[gray!75,line width=0.5pt] plot[smooth] coordinates {(1.448,1.052) (1.399,1.074) (1.347,1.098) (1.290,1.123) (1.229,1.150) (1.162,1.178) (1.088,1.209) (1.007,1.241) (0.918,1.274) (0.819,1.307) (0.710,1.340) (0.588,1.371) (0.454,1.397)};
\draw[gray!75,line width=0.5pt] plot[smooth] coordinates {(-3.913,-2.843) (-3.781,-2.902) (-3.639,-2.965) (-3.486,-3.032) (-3.320,-3.105) (-3.140,-3.182) (-2.941,-3.264) (-2.723,-3.349) (-2.482,-3.438) (-2.214,-3.528) (-1.918,-3.617) (-1.588,-3.699) (-1.225,-3.770)};
\draw[gray!75,line width=0.5pt] plot[smooth] coordinates {(-3.913,-2.843) (-3.694,-3.238) (-3.474,-3.589) (-3.259,-3.900) (-3.054,-4.176) (-2.860,-4.423) (-2.678,-4.643) (-2.508,-4.843) (-2.351,-5.023) (-2.204,-5.189) (-2.068,-5.341) (-1.941,-5.481) (-1.824,-5.613)};
\draw[gray!75,line width=0.5pt] plot[smooth] coordinates {(0.978,0.711) (0.924,0.808) (0.869,0.894) (0.816,0.971) (0.764,1.040) (0.716,1.101) (0.670,1.156) (0.627,1.205) (0.588,1.250) (0.550,1.292) (0.516,1.329) (0.484,1.364) (0.454,1.397)};
\draw[gray!75,line width=0.5pt] plot[smooth] coordinates {(0.978,0.711) (0.945,0.727) (0.909,0.745) (0.870,0.764) (0.828,0.784) (0.783,0.805) (0.733,0.827) (0.678,0.849) (0.618,0.873) (0.552,0.896) (0.479,0.919) (0.399,0.940) (0.311,0.959)};
\draw[gray!75,line width=0.5pt] plot[smooth] coordinates {(-5.825,-4.232) (-5.628,-4.319) (-5.417,-4.413) (-5.190,-4.514) (-4.943,-4.622) (-4.673,-4.736) (-4.378,-4.858) (-4.054,-4.986) (-3.694,-5.118) (-3.296,-5.252) (-2.854,-5.384) (-2.364,-5.506) (-1.824,-5.613)};
\draw[gray!75,line width=0.5pt] plot[smooth] coordinates {(-5.825,-4.232) (-5.361,-4.629) (-4.833,-5.032) (-4.238,-5.431) (-3.575,-5.809) (-2.854,-6.148) (-2.090,-6.433) (-1.304,-6.651) (-0.522,-6.801) (0.236,-6.884) (0.952,-6.912) (1.616,-6.896) (2.225,-6.848)};
\draw[gray!75,line width=0.5pt] plot[smooth] coordinates {(-5.825,-4.232) (-6.059,-3.668) (-6.279,-3.041) (-6.474,-2.352) (-6.629,-1.605) (-6.729,-0.815) (-6.764,0.000) (-6.729,0.815) (-6.629,1.605) (-6.474,2.352) (-6.279,3.041) (-6.059,3.668) (-5.825,4.232)};
\draw[gray!75,line width=0.5pt] plot[smooth] coordinates {(2.225,6.848) (2.015,6.803) (1.789,6.755) (1.545,6.702) (1.282,6.644) (0.997,6.579) (0.687,6.504) (0.349,6.416) (-0.019,6.312) (-0.420,6.186) (-0.855,6.033) (-1.324,5.844) (-1.824,5.613)};
\draw[gray!75,line width=0.5pt] plot[smooth] coordinates {(2.225,6.848) (1.616,6.896) (0.952,6.912) (0.236,6.884) (-0.522,6.801) (-1.304,6.651) (-2.090,6.433) (-2.854,6.148) (-3.575,5.809) (-4.238,5.431) (-4.833,5.032) (-5.361,4.629) (-5.825,4.232)};
\draw[gray!75,line width=0.5pt] plot[smooth] coordinates {(-0.374,-1.150) (-0.337,-1.144) (-0.297,-1.137) (-0.255,-1.129) (-0.209,-1.121) (-0.160,-1.111) (-0.107,-1.099) (-0.049,-1.086) (0.013,-1.069) (0.080,-1.049) (0.153,-1.025) (0.230,-0.995) (0.311,-0.959)};
\draw[gray!75,line width=0.5pt] plot[smooth] coordinates {(-0.374,-1.150) (-0.273,-1.197) (-0.177,-1.235) (-0.089,-1.265) (-0.007,-1.291) (0.068,-1.312) (0.137,-1.329) (0.201,-1.344) (0.260,-1.357) (0.314,-1.368) (0.364,-1.379) (0.411,-1.388) (0.454,-1.397)};
\draw[gray!75,line width=0.5pt] plot[smooth] coordinates {(1.495,4.600) (1.085,4.790) (0.701,4.945) (0.344,5.071) (0.016,5.174) (-0.286,5.259) (-0.563,5.331) (-0.817,5.392) (-1.051,5.446) (-1.267,5.493) (-1.466,5.536) (-1.651,5.576) (-1.824,5.613)};
\draw[gray!75,line width=0.5pt] plot[smooth] coordinates {(1.495,4.600) (1.353,4.570) (1.202,4.538) (1.038,4.502) (0.861,4.463) (0.670,4.420) (0.461,4.369) (0.234,4.310) (-0.013,4.240) (-0.282,4.156) (-0.574,4.053) (-0.889,3.926) (-1.225,3.770)};
\draw[gray!75,line width=0.5pt] plot[smooth] coordinates {(-0.553,-1.702) (-0.501,-1.691) (-0.444,-1.679) (-0.384,-1.667) (-0.318,-1.652) (-0.247,-1.636) (-0.170,-1.618) (-0.086,-1.596) (0.006,-1.570) (0.105,-1.539) (0.214,-1.501) (0.330,-1.454) (0.454,-1.397)};
\draw[gray!75,line width=0.5pt] plot[smooth] coordinates {(-0.553,-1.702) (-0.402,-1.773) (-0.260,-1.830) (-0.128,-1.876) (-0.006,-1.914) (0.105,-1.946) (0.208,-1.972) (0.302,-1.995) (0.389,-2.015) (0.468,-2.032) (0.542,-2.048) (0.611,-2.063) (0.675,-2.076)};
\draw[gray!75,line width=0.5pt] plot[smooth] coordinates {(1.004,3.090) (0.729,3.218) (0.471,3.322) (0.231,3.406) (0.011,3.476) (-0.192,3.533) (-0.378,3.581) (-0.549,3.622) (-0.706,3.658) (-0.851,3.690) (-0.985,3.719) (-1.109,3.746) (-1.225,3.770)};
\draw[gray!75,line width=0.5pt] plot[smooth] coordinates {(1.004,3.090) (0.909,3.070) (0.807,3.048) (0.697,3.025) (0.579,2.999) (0.450,2.969) (0.310,2.935) (0.157,2.896) (-0.009,2.849) (-0.190,2.792) (-0.386,2.723) (-0.598,2.638) (-0.823,2.533)};
\draw[gray!75,line width=0.5pt] plot[smooth] coordinates {(-0.823,-2.533) (-0.745,-2.516) (-0.662,-2.499) (-0.572,-2.479) (-0.474,-2.458) (-0.369,-2.434) (-0.254,-2.406) (-0.129,-2.374) (0.007,-2.335) (0.156,-2.289) (0.316,-2.232) (0.490,-2.162) (0.675,-2.076)};
\draw[gray!75,line width=0.5pt] plot[smooth] coordinates {(-0.823,-2.533) (-0.598,-2.638) (-0.386,-2.723) (-0.190,-2.792) (-0.009,-2.849) (0.157,-2.896) (0.310,-2.935) (0.450,-2.969) (0.579,-2.999) (0.697,-3.025) (0.807,-3.048) (0.909,-3.070) (1.004,-3.090)};
\draw[gray!75,line width=0.5pt] plot[smooth] coordinates {(0.675,2.076) (0.490,2.162) (0.316,2.232) (0.156,2.289) (0.007,2.335) (-0.129,2.374) (-0.254,2.406) (-0.369,2.434) (-0.474,2.458) (-0.572,2.479) (-0.662,2.499) (-0.745,2.516) (-0.823,2.533)};
\draw[gray!75,line width=0.5pt] plot[smooth] coordinates {(0.675,2.076) (0.611,2.063) (0.542,2.048) (0.468,2.032) (0.389,2.015) (0.302,1.995) (0.208,1.972) (0.105,1.946) (-0.006,1.914) (-0.128,1.876) (-0.260,1.830) (-0.402,1.773) (-0.553,1.702)};
\draw[gray!75,line width=0.5pt] plot[smooth] coordinates {(-1.225,-3.770) (-1.109,-3.746) (-0.985,-3.719) (-0.851,-3.690) (-0.706,-3.658) (-0.549,-3.622) (-0.378,-3.581) (-0.192,-3.533) (0.011,-3.476) (0.231,-3.406) (0.471,-3.322) (0.729,-3.218) (1.004,-3.090)};
\draw[gray!75,line width=0.5pt] plot[smooth] coordinates {(-1.225,-3.770) (-0.889,-3.926) (-0.574,-4.053) (-0.282,-4.156) (-0.013,-4.240) (0.234,-4.310) (0.461,-4.369) (0.670,-4.420) (0.861,-4.463) (1.038,-4.502) (1.202,-4.538) (1.353,-4.570) (1.495,-4.600)};
\draw[gray!75,line width=0.5pt] plot[smooth] coordinates {(0.454,1.397) (0.330,1.454) (0.214,1.501) (0.105,1.539) (0.006,1.570) (-0.086,1.596) (-0.170,1.618) (-0.247,1.636) (-0.318,1.652) (-0.384,1.667) (-0.444,1.679) (-0.501,1.691) (-0.553,1.702)};
\draw[gray!75,line width=0.5pt] plot[smooth] coordinates {(0.454,1.397) (0.411,1.388) (0.364,1.379) (0.314,1.368) (0.260,1.357) (0.201,1.344) (0.137,1.329) (0.068,1.312) (-0.007,1.291) (-0.089,1.265) (-0.177,1.235) (-0.273,1.197) (-0.374,1.150)};
\draw[gray!75,line width=0.5pt] plot[smooth] coordinates {(-1.824,-5.613) (-1.651,-5.576) (-1.466,-5.536) (-1.267,-5.493) (-1.051,-5.446) (-0.817,-5.392) (-0.563,-5.331) (-0.286,-5.259) (0.016,-5.174) (0.344,-5.071) (0.701,-4.945) (1.085,-4.790) (1.495,-4.600)};
\draw[gray!75,line width=0.5pt] plot[smooth] coordinates {(-1.824,-5.613) (-1.324,-5.844) (-0.855,-6.033) (-0.420,-6.186) (-0.019,-6.312) (0.349,-6.416) (0.687,-6.504) (0.997,-6.579) (1.282,-6.644) (1.545,-6.702) (1.789,-6.755) (2.015,-6.803) (2.225,-6.848)};
\draw[gray!75,line width=0.5pt] plot[smooth] coordinates {(0.311,0.959) (0.230,0.995) (0.153,1.025) (0.080,1.049) (0.013,1.069) (-0.049,1.086) (-0.107,1.099) (-0.160,1.111) (-0.209,1.121) (-0.255,1.129) (-0.297,1.137) (-0.337,1.144) (-0.374,1.150)};
\draw[gray!75,line width=0.5pt] plot[smooth] coordinates {(0.311,0.959) (0.304,0.934) (0.295,0.908) (0.286,0.881) (0.276,0.851) (0.266,0.818) (0.254,0.781) (0.241,0.740) (0.225,0.693) (0.207,0.637) (0.184,0.565) (0.150,0.461) (-0.000,0.000)};
\draw[gray!75,line width=0.5pt] plot[smooth] coordinates {(0.311,-0.959) (0.399,-0.940) (0.479,-0.919) (0.552,-0.896) (0.618,-0.873) (0.678,-0.849) (0.733,-0.827) (0.783,-0.805) (0.828,-0.784) (0.870,-0.764) (0.909,-0.745) (0.945,-0.727) (0.978,-0.711)};
\draw[gray!75,line width=0.5pt] plot[smooth] coordinates {(0.311,-0.959) (0.304,-0.934) (0.295,-0.908) (0.286,-0.881) (0.276,-0.851) (0.266,-0.818) (0.254,-0.781) (0.241,-0.740) (0.225,-0.693) (0.207,-0.637) (0.184,-0.565) (0.150,-0.461) (-0.000,0.000)};
\draw[gray!75,line width=0.5pt] plot[smooth] coordinates {(-1.824,5.613) (-2.364,5.506) (-2.854,5.384) (-3.296,5.252) (-3.694,5.118) (-4.054,4.986) (-4.378,4.858) (-4.673,4.736) (-4.943,4.622) (-5.190,4.514) (-5.417,4.413) (-5.628,4.319) (-5.825,4.232)};
\draw[gray!75,line width=0.5pt] plot[smooth] coordinates {(-1.824,5.613) (-1.941,5.481) (-2.068,5.341) (-2.204,5.189) (-2.351,5.023) (-2.508,4.843) (-2.678,4.643) (-2.860,4.423) (-3.054,4.176) (-3.259,3.900) (-3.474,3.589) (-3.694,3.238) (-3.913,2.843)};
\draw[gray!75,line width=0.5pt] plot[smooth] coordinates {(0.454,-1.397) (0.484,-1.364) (0.516,-1.329) (0.550,-1.292) (0.588,-1.250) (0.627,-1.205) (0.670,-1.156) (0.716,-1.101) (0.764,-1.040) (0.816,-0.971) (0.869,-0.894) (0.924,-0.808) (0.978,-0.711)};
\draw[gray!75,line width=0.5pt] plot[smooth] coordinates {(0.454,-1.397) (0.588,-1.371) (0.710,-1.340) (0.819,-1.307) (0.918,-1.274) (1.007,-1.241) (1.088,-1.209) (1.162,-1.178) (1.229,-1.150) (1.290,-1.123) (1.347,-1.098) (1.399,-1.074) (1.448,-1.052)};
\draw[gray!75,line width=0.5pt] plot[smooth] coordinates {(-1.225,3.770) (-1.588,3.699) (-1.918,3.617) (-2.214,3.528) (-2.482,3.438) (-2.723,3.349) (-2.941,3.264) (-3.140,3.182) (-3.320,3.105) (-3.486,3.032) (-3.639,2.965) (-3.781,2.902) (-3.913,2.843)};
\draw[gray!75,line width=0.5pt] plot[smooth] coordinates {(-1.225,3.770) (-1.304,3.682) (-1.389,3.588) (-1.481,3.486) (-1.579,3.375) (-1.685,3.253) (-1.799,3.119) (-1.921,2.971) (-2.052,2.806) (-2.189,2.620) (-2.334,2.411) (-2.481,2.175) (-2.629,1.910)};
\draw[gray!75,line width=0.5pt] plot[smooth] coordinates {(0.675,-2.076) (0.718,-2.028) (0.765,-1.976) (0.816,-1.920) (0.870,-1.858) (0.928,-1.792) (0.991,-1.718) (1.058,-1.636) (1.130,-1.545) (1.206,-1.443) (1.286,-1.328) (1.367,-1.198) (1.448,-1.052)};
\draw[gray!75,line width=0.5pt] plot[smooth] coordinates {(0.675,-2.076) (0.875,-2.037) (1.056,-1.992) (1.219,-1.943) (1.367,-1.893) (1.499,-1.844) (1.620,-1.797) (1.729,-1.752) (1.828,-1.710) (1.920,-1.670) (2.004,-1.633) (2.082,-1.598) (2.155,-1.566)};
\draw[gray!75,line width=0.5pt] plot[smooth] coordinates {(-0.823,2.533) (-1.067,2.485) (-1.288,2.430) (-1.488,2.370) (-1.667,2.310) (-1.829,2.250) (-1.976,2.193) (-2.109,2.138) (-2.231,2.086) (-2.342,2.037) (-2.445,1.992) (-2.540,1.949) (-2.629,1.910)};
\draw[gray!75,line width=0.5pt] plot[smooth] coordinates {(-0.823,2.533) (-0.876,2.474) (-0.933,2.410) (-0.995,2.342) (-1.061,2.267) (-1.132,2.186) (-1.209,2.096) (-1.291,1.996) (-1.378,1.885) (-1.471,1.760) (-1.568,1.620) (-1.667,1.461) (-1.766,1.283)};
\draw[gray!75,line width=0.5pt] plot[smooth] coordinates {(1.004,-3.090) (1.069,-3.018) (1.139,-2.941) (1.214,-2.857) (1.294,-2.766) (1.381,-2.666) (1.475,-2.557) (1.575,-2.435) (1.682,-2.300) (1.795,-2.147) (1.913,-1.976) (2.034,-1.783) (2.155,-1.566)};
\draw[gray!75,line width=0.5pt] plot[smooth] coordinates {(1.004,-3.090) (1.302,-3.032) (1.572,-2.964) (1.815,-2.892) (2.034,-2.818) (2.232,-2.745) (2.411,-2.675) (2.573,-2.608) (2.722,-2.545) (2.858,-2.485) (2.983,-2.430) (3.099,-2.378) (3.207,-2.330)};
\draw[gray!75,line width=0.5pt] plot[smooth] coordinates {(-0.553,1.702) (-0.717,1.670) (-0.865,1.633) (-0.999,1.593) (-1.120,1.552) (-1.229,1.512) (-1.327,1.473) (-1.417,1.436) (-1.499,1.402) (-1.574,1.369) (-1.643,1.338) (-1.707,1.310) (-1.766,1.283)};
\draw[gray!75,line width=0.5pt] plot[smooth] coordinates {(-0.553,1.702) (-0.589,1.663) (-0.628,1.620) (-0.669,1.574) (-0.714,1.524) (-0.762,1.469) (-0.813,1.408) (-0.869,1.342) (-0.928,1.267) (-0.990,1.183) (-1.055,1.089) (-1.122,0.983) (-1.188,0.863)};
\draw[gray!75,line width=0.5pt] plot[smooth] coordinates {(1.495,-4.600) (1.591,-4.493) (1.695,-4.377) (1.807,-4.253) (1.927,-4.117) (2.056,-3.969) (2.195,-3.806) (2.344,-3.625) (2.503,-3.423) (2.671,-3.196) (2.847,-2.941) (3.027,-2.654) (3.207,-2.330)};
\draw[gray!75,line width=0.5pt] plot[smooth] coordinates {(1.495,-4.600) (1.938,-4.513) (2.339,-4.412) (2.702,-4.305) (3.028,-4.195) (3.322,-4.086) (3.589,-3.982) (3.830,-3.882) (4.051,-3.788) (4.254,-3.700) (4.440,-3.617) (4.613,-3.540) (4.774,-3.469)};
\draw[gray!75,line width=0.5pt] plot[smooth] coordinates {(-0.374,1.150) (-0.483,1.128) (-0.582,1.103) (-0.672,1.076) (-0.753,1.048) (-0.826,1.021) (-0.892,0.994) (-0.953,0.969) (-1.008,0.945) (-1.058,0.923) (-1.105,0.902) (-1.148,0.882) (-1.188,0.863)};
\draw[gray!75,line width=0.5pt] plot[smooth] coordinates {(-0.374,1.150) (-0.400,1.123) (-0.428,1.095) (-0.457,1.064) (-0.489,1.030) (-0.523,0.993) (-0.560,0.952) (-0.598,0.907) (-0.639,0.857) (-0.682,0.802) (-0.726,0.739) (-0.771,0.670) (-0.815,0.592)};
\draw[gray!75,line width=0.5pt] plot[smooth] coordinates {(2.225,-6.848) (2.369,-6.688) (2.523,-6.516) (2.689,-6.330) (2.868,-6.129) (3.060,-5.908) (3.267,-5.665) (3.489,-5.396) (3.726,-5.095) (3.976,-4.758) (4.238,-4.378) (4.506,-3.950) (4.774,-3.469)};
\draw[gray!75,line width=0.5pt] plot[smooth] coordinates {(-0.815,0.592) (-0.795,0.577) (-0.773,0.561) (-0.749,0.544) (-0.724,0.526) (-0.696,0.505) (-0.665,0.483) (-0.630,0.458) (-0.590,0.428) (-0.542,0.394) (-0.481,0.349) (-0.392,0.285) (-0.000,0.000)};
\node[o] at (7.200,0.000) {}; \node[font=\scriptsize,red!80!black] at (7.660,0.300) {$(0,0)$};
\node[v] at (-1.209,0.000) {};
\node[font=\tiny,inner sep=0.5pt,  text opacity=1] at (-0.989,0.170) {$\scriptscriptstyle(0,1)$};
\node[v] at (4.837,0.000) {};
\node[font=\tiny,inner sep=0.5pt,  text opacity=1] at (5.057,0.170) {$\scriptscriptstyle(0,2)$};
\node[v] at (-1.790,0.000) {};
\node[font=\tiny,inner sep=0.5pt,  text opacity=1] at (-1.570,0.170) {$\scriptscriptstyle(0,3)$};
\node[v] at (3.249,0.000) {};
\node[font=\tiny,inner sep=0.5pt,  text opacity=1] at (3.469,0.170) {$\scriptscriptstyle(0,4)$};
\node[v] at (-2.663,0.000) {};
\node[font=\tiny,inner sep=0.5pt,  text opacity=1] at (-2.443,0.170) {$\scriptscriptstyle(0,5)$};
\node[v] at (2.183,0.000) {};
\node[font=\tiny,inner sep=0.5pt,  text opacity=1] at (2.403,0.170) {$\scriptscriptstyle(0,6)$};
\node[v] at (-3.965,0.000) {};
\node[font=\tiny,inner sep=0.5pt,  text opacity=1] at (-3.745,0.170) {$\scriptscriptstyle(0,7)$};
\node[v] at (1.469,0.000) {};
\node[font=\tiny,inner sep=0.5pt,  text opacity=1] at (1.689,0.170) {$\scriptscriptstyle(0,8)$};
\node[v] at (-5.901,0.000) {};
\node[font=\tiny,inner sep=0.5pt,  text opacity=1] at (-5.681,0.170) {$\scriptscriptstyle(0,9)$};
\node[v] at (1.008,0.000) {};
\node[font=\tiny,inner sep=0.5pt,  text opacity=1] at (0.928,0.170) {$\scriptscriptstyle(0,10)$};
\node[v] at (-0.815,-0.592) {};
\node[font=\tiny,inner sep=0.5pt,  text opacity=1] at (-0.595,-0.422) {$\scriptscriptstyle(1,0)$};
\node[v] at (4.774,3.469) {};
\node[font=\tiny,inner sep=0.5pt,  text opacity=1] at (4.994,3.639) {$\scriptscriptstyle(1,1)$};
\node[v] at (-1.188,-0.863) {};
\node[font=\tiny,inner sep=0.5pt,  text opacity=1] at (-1.568,-0.693) {$\scriptscriptstyle(1,2)$};
\node[v] at (3.207,2.330) {};
\node[font=\tiny,inner sep=0.5pt,  text opacity=1] at (3.427,2.500) {$\scriptscriptstyle(1,3)$};
\node[v] at (-1.766,-1.283) {};
\node[font=\tiny,inner sep=0.5pt,  text opacity=1] at (-1.546,-1.113) {$\scriptscriptstyle(1,4)$};
\node[v] at (2.155,1.566) {};
\node[font=\tiny,inner sep=0.5pt,  text opacity=1] at (2.375,1.736) {$\scriptscriptstyle(1,5)$};
\node[v] at (-2.629,-1.910) {};
\node[font=\tiny,inner sep=0.5pt,  text opacity=1] at (-2.409,-1.740) {$\scriptscriptstyle(1,6)$};
\node[v] at (1.448,1.052) {};
\node[font=\tiny,inner sep=0.5pt,  text opacity=1] at (1.668,1.222) {$\scriptscriptstyle(1,7)$};
\node[v] at (-3.913,-2.843) {};
\node[font=\tiny,inner sep=0.5pt,  text opacity=1] at (-3.693,-2.673) {$\scriptscriptstyle(1,8)$};
\node[v] at (0.978,0.711) {};
\node[font=\tiny,inner sep=0.5pt,  text opacity=1] at (1.198,0.881) {$\scriptscriptstyle(1,9)$};
\node[o] at (-5.825,-4.232) {}; \node[font=\scriptsize,red!80!black] at (-6.116,-4.544) {$(1,10)$};
\node[o] at (2.225,6.848) {}; \node[font=\scriptsize,red!80!black] at (2.336,7.190) {$(2,0)$};
\node[v] at (-0.374,-1.150) {};
\node[font=\tiny,inner sep=0.5pt,  text opacity=1] at (-0.154,-0.980) {$\scriptscriptstyle(2,1)$};
\node[v] at (1.495,4.600) {};
\node[font=\tiny,inner sep=0.5pt,  text opacity=1] at (1.715,4.770) {$\scriptscriptstyle(2,2)$};
\node[v] at (-0.553,-1.702) {};
\node[font=\tiny,inner sep=0.5pt,  text opacity=1] at (-0.333,-1.532) {$\scriptscriptstyle(2,3)$};
\node[v] at (1.004,3.090) {};
\node[font=\tiny,inner sep=0.5pt,  text opacity=1] at (1.224,3.260) {$\scriptscriptstyle(2,4)$};
\node[v] at (-0.823,-2.533) {};
\node[font=\tiny,inner sep=0.5pt,  text opacity=1] at (-0.603,-2.363) {$\scriptscriptstyle(2,5)$};
\node[v] at (0.675,2.076) {};
\node[font=\tiny,inner sep=0.5pt,  text opacity=1] at (0.895,2.246) {$\scriptscriptstyle(2,6)$};
\node[v] at (-1.225,-3.770) {};
\node[font=\tiny,inner sep=0.5pt,  text opacity=1] at (-1.005,-3.600) {$\scriptscriptstyle(2,7)$};
\node[v] at (0.454,1.397) {};
\node[font=\tiny,inner sep=0.5pt,  text opacity=1] at (0.674,1.567) {$\scriptscriptstyle(2,8)$};
\node[v] at (-1.824,-5.613) {};
\node[font=\tiny,inner sep=0.5pt,  text opacity=1] at (-1.604,-5.443) {$\scriptscriptstyle(2,9)$};
\node[v] at (0.311,0.959) {};
\node[font=\tiny,inner sep=0.5pt,  text opacity=1] at (0.531,1.129) {$\scriptscriptstyle(2,10)$};
\node[v] at (0.311,-0.959) {};
\node[font=\tiny,inner sep=0.5pt,  text opacity=1] at (0.531,-0.789) {$\scriptscriptstyle(3,0)$};
\node[v] at (-1.824,5.613) {};
\node[font=\tiny,inner sep=0.5pt,  text opacity=1] at (-1.604,5.783) {$\scriptscriptstyle(3,1)$};
\node[v] at (0.454,-1.397) {};
\node[font=\tiny,inner sep=0.5pt,  text opacity=1] at (0.674,-1.227) {$\scriptscriptstyle(3,2)$};
\node[v] at (-1.225,3.770) {};
\node[font=\tiny,inner sep=0.5pt,  text opacity=1] at (-1.005,3.940) {$\scriptscriptstyle(3,3)$};
\node[v] at (0.675,-2.076) {};
\node[font=\tiny,inner sep=0.5pt,  text opacity=1] at (0.895,-1.906) {$\scriptscriptstyle(3,4)$};
\node[v] at (-0.823,2.533) {};
\node[font=\tiny,inner sep=0.5pt,  text opacity=1] at (-0.603,2.703) {$\scriptscriptstyle(3,5)$};
\node[v] at (1.004,-3.090) {};
\node[font=\tiny,inner sep=0.5pt,  text opacity=1] at (1.224,-2.920) {$\scriptscriptstyle(3,6)$};
\node[v] at (-0.553,1.702) {};
\node[font=\tiny,inner sep=0.5pt,  text opacity=1] at (-0.333,1.872) {$\scriptscriptstyle(3,7)$};
\node[v] at (1.495,-4.600) {};
\node[font=\tiny,inner sep=0.5pt,  text opacity=1] at (1.715,-4.430) {$\scriptscriptstyle(3,8)$};
\node[v] at (-0.374,1.150) {};
\node[font=\tiny,inner sep=0.5pt,  text opacity=1] at (-0.154,1.320) {$\scriptscriptstyle(3,9)$};
\node[o] at (2.225,-6.848) {}; \node[font=\scriptsize,red!80!black] at (2.336,-7.190) {$(3,10)$};
\node[o] at (-5.825,4.232) {}; \node[font=\scriptsize,red!80!black] at (-6.116,4.544) {$(4,0)$};
\node[v] at (0.978,-0.711) {};
\node[font=\tiny,inner sep=0.5pt,  text opacity=1] at (1.198,-0.541) {$\scriptscriptstyle(4,1)$};
\node[v] at (-3.913,2.843) {};
\node[font=\tiny,inner sep=0.5pt,  text opacity=1] at (-3.693,3.013) {$\scriptscriptstyle(4,2)$};
\node[v] at (1.448,-1.052) {};
\node[font=\tiny,inner sep=0.5pt,  text opacity=1] at (1.668,-0.882) {$\scriptscriptstyle(4,3)$};
\node[v] at (-2.629,1.910) {};
\node[font=\tiny,inner sep=0.5pt,  text opacity=1] at (-2.409,2.080) {$\scriptscriptstyle(4,4)$};
\node[v] at (2.155,-1.566) {};
\node[font=\tiny,inner sep=0.5pt,  text opacity=1] at (2.375,-1.396) {$\scriptscriptstyle(4,5)$};
\node[v] at (-1.766,1.283) {};
\node[font=\tiny,inner sep=0.5pt,  text opacity=1] at (-1.546,1.453) {$\scriptscriptstyle(4,6)$};
\node[v] at (3.207,-2.330) {};
\node[font=\tiny,inner sep=0.5pt,  text opacity=1] at (3.427,-2.160) {$\scriptscriptstyle(4,7)$};
\node[v] at (-1.188,0.863) {};
\node[font=\tiny,inner sep=0.5pt,  text opacity=1] at (-0.968,1.033) {$\scriptscriptstyle(4,8)$};
\node[v] at (4.774,-3.469) {};
\node[font=\tiny,inner sep=0.5pt,  text opacity=1] at (4.994,-3.299) {$\scriptscriptstyle(4,9)$};
\node[v] at (-0.815,0.592) {};
\node[font=\tiny,inner sep=0.5pt,  text opacity=1] at (-0.595,0.762) {$\scriptscriptstyle(4,10)$};
\node[font=\tiny,inner sep=1.5pt,fill=violet,fill opacity=0.7,text opacity=1] at (-0.000,0.000) {}; 
\node[above right=1pt,font=\footnotesize, violet] at (-0.000,0.000) {$a$};
\end{tikzpicture}
    \caption{1-Skeleton of the Reuleaux polyhedron whose diameter graph is $Q_{5,\, 10}$.}
    \label{Fig:Q510}
\end{figure}
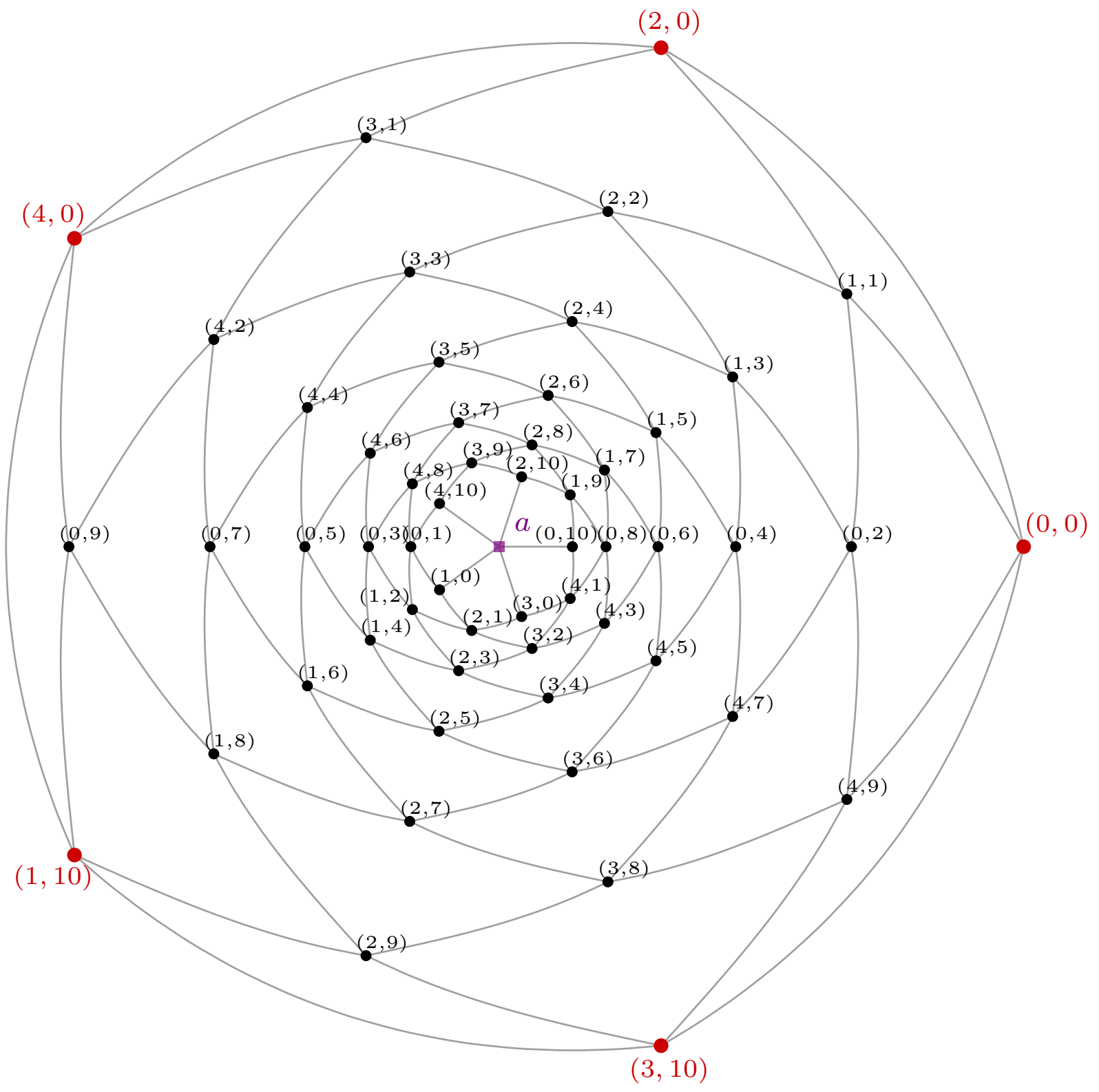

Nevertheless, we would like to acknowledge that the family of twisted prismoids was discovered thanks to Damásdi’s embedding of diameter graphs on the projective plane, and since we do believe in the correctness of Theorem~\ref{Thm:Damasdi}, which we think is currently undergoing peer review and is expected to be revised before the publication of this paper, we have decided to keep the present approach.

\subsection{On Brass Conjecture}\label{Sec:BrassConjecture}
In this subsection we prove Theorem \ref{Thm:Brass}, confirming Brass's conjecture.

\bigskip
\begin{lemma}\label{Lemma:unique}
    Let $V\subset \mathbb{R}^3$ be a finite point set with unit diameter and Borsuk number 4. Then, there exist a unique subset $V_1\subset V$ with unit diameter such that $\Bor (V_1)$ is 4-critical.
\end{lemma} 
\begin{proof}
    Since $\chi(\Bor (V))=\bn(V)=4$, Proposition \ref{Prop:critical} guarantees the existence of a subset $V_1\subset V$ with unit diameter such that $\Bor (V_1)$ is 4-critical. It remains to prove that $V_1$ is unique. We proceed by contradiction. Suppose there is a subset $V_2\subset V$ with unit diameter, distinct from $V_1$, such that $\Bor (V_2)$ is 4-critical. 
    
    Since $V_1\cup V_2\subseteq V$, the set $V_1\cup V_2$ has diameter 1, which implies $V_2\subset \mathcal{B}(V_1)$ and $V_1\subset \mathcal{B}(V_2)$ by definition. Furthermore, by Theorem \ref{finiteC} ((ii)$\implies$(iii)), both $\mathcal{B}(V_1)$ and $\mathcal{B}(V_2)$ are Reuleaux polyhedra. 

Since both $\Bor (V_1)$ and $\Bor (V_2)$ are 4-critical, $V_1\backslash V_2\neq \emptyset$. Let $v_1\in V_1\backslash V_2$, and let $\epsilon$ be half of the distance from $v_1$ to $V_2$. By Theorem \ref{MeissnerThm2025},   there is a cover $P_{\epsilon, v_1}=\{P_1,  P_2,  P_3,  P_4\}$ of $\mathcal{B}(V_1)$, with $\diam(P_i)<1$ for $i=2,3,4$. Clearly $V_2 \subset P_2\cup P_3\cup P_4$, then $\chi(\Bor(V_2))=\bn (V_2)\leq 3$, a contradiction.
\end{proof}

\begin{proof}[Proof of Theorem \ref{Thm:Brass}]
       Let $V\subset \mathbb{R}^3$ be a unit diameter finite point set. If $\chi (\Bor (V))\leq 3$, we are done. Suppose that $\chi (\Bor (V))= 4$ ($\chi (\Bor (V))\leq 4$ since Borsuk's conjecture holds in $\mathbb{R}^3$). By Lemma \ref{Lemma:unique}, there exist a unique unit diameter set subset $V_1\subset V$ such that $\Bor (V_1)$ is 4-critical. 
       
To conclude the proof, it is enough to prove the following claim.
       
\textbf{Claim:} Let $x\in V$, then $\chi (\Bor (V)-\{x\})=3$ if and only if $x\in V_1$.

         \emph{Necessity:} By contradiction, if $x\notin V_1$, then $V_1\subset V\setminus \{x\}$, which implies that $\chi (\Bor (V)-\{x\})=4$, a contradiction.
       
       \emph{Sufficiency:} By contradiction, suppose there exist $x\in V_1$ with $\chi (\Bor (V)-\{x\})\neq 3$, then $\chi (\Bor (\Bor (V)-\{x\}))=4$, because deleting a vertex decreases the chromatic number by at most one. Therefore, by Proposition \ref{Prop:critical}, there is a 4-critical subgraph of $\Bor (V)-\{x\}$ with set of vertices different than $V_1$, a contradiction with the uniqueness of $V_1$.
       
Hence removing any point of $V_1$ leaves a $3$-colorable graph.
\end{proof}
\section{Brass-type problems}\label{Sec:Brasstype}

In this section, we address some cases of Question~\ref{Question:Main} and establish some general properties. We note that, as far as we know, Question~\ref{Question:Main} has only been posed for $k=2,3$ and $d=3$ (as mentioned in the Introduction), but not in this general setting. Therefore, we refer to questions of this type as \emph{Brass-type problems}.

 Let $S\subset\mathbb{R}^d$ be a finite point set. Denote by $\Br_S(k)$ the minimum number of points that must be removed from $S$, such that the subgraph of the diameter graph induced by the remaining points is $k$-colorable.  Let $\Phi_d$ be the family of finite subsets of $\mathbb{R}^d$. We define the \emph{Brass number} $\Br(k,d)$ as:
\[
\Br(k,d):=\sup_{S\in \Phi_d} \Br_S(k). 
\]

Thus, we can rephrase Question \ref{Question:Main} as follows:\\
\begin{question}\label{Question:main}
        Let $k$ and $d$ be fixed positive integers. When is $\Br(k,d)$ finite?
\end{question}

We present some straightforward properties of $\Br(k,d)$.\\

\smallskip
\begin{proposition}\label{Prop:dim}
    Let $k$ and $d$ be fixed positive integers. Then, $\Br(k,d)\leq \Br(k,j)$ for all $j>d$. In particular, if $\Br(k,d)$ is not finite, then neither is $\Br(k,j)$.
\end{proposition}
\begin{proof}
    Let $\phi_j:\mathbb{R}^d\rightarrow \mathbb{R}^j$ be an isometric embedding. For every finite set $S\subset \mathbb{R}^d$, it is easy to see that $\Bor(\phi_j(S))\cong\Bor(S)$, then $\Br_S(k)=\Br_{\phi(S)}(k)$ for every positive integer $k$, which, by definition, implies \[\Br(k, d)\leq \Br(k, j).\qedhere\] 
\end{proof}
\smallskip
\begin{proposition}\label{Prop:color}
    Let $k$ and $d$ be fixed positive integers. Then $\Br(k,d)\geq \Br(i,d)$ for all $i>k$. In particular, if $\Br(k,d)$ is finite, then so is $\Br(i,d)$. 
\end{proposition}
\begin{proof}
    Since for every finite set $S\subset \mathbb{R}^d$, there exists a subset $S_1\subset S$, with $|S_1|=\Br_S(k)$, such that the subgraph $G_1$ of $\Bor(S)$ induced by $S\setminus S_1$ is $k$-colorable, then $G_1$ is particularly $i$-colorable, for every $i>k$. 

    Therefore, $\Br_S(k)\geq \Br_S(i)$, which, by definition, implies
    \[\Br(k,d)\geq \Br(i,d).\qedhere\] 
\end{proof}
\bigskip

\begin{proposition}\label{Prop-simplex}
Let $d$ be a fixed positive integer number. Then, $\Br(d,d)>0$.
\end{proposition}
\begin{proof}
Let $S\subset \mathbb{R}^d$ be the vertex set of a regular $d-$simplex. Since $S$ has diameter graph isomorphic to the $(d+1)$-clique (the complete graph of $d$ vertices), then $\chi(\Bor(S))=d+1$. The latter implies that exactly one point must be removed from $S$ to obtain a $d$-colorable diameter graph, that is $\Br_S(d)=1$. Therefore,
\[\Br(d,d)\geq\Br_S(d)=1\qedhere\]
\end{proof}
\bigskip

\begin{proposition}\label{Prop:1-color}
    Let $d$ be a fixed positive integer number. Then, $\Br(1,d)=1$ if and only if $d=1$, otherwise $\Br(1,d)$ is not finite.
\end{proposition}
\begin{proof}
    Suppose $d=1$. Since for every finite point set $S\subset \mathbb{R}$ with at least two points, $\Bor(S)$ only has one edge, joining the two extremal points of $S$, then $\chi(\Bor(S))=2$. Thus, by removing only one of these extremal points, the remaining graph would have no edges, and hence it would be $1$-colorable. If $|S|=1$, clearly $\Br_S(1)=0$. Therefore, $\Br(1,1)=1$.

    Now suppose $\Br(1,d)=1$ for some $d$. Let $r$ be a positive integer, and let $S_r\subset \mathbb{R}^2$ be the set of vertices of a regular polygon with $2(r+1)$ vertices.
    
    Clearly, after removing any $r$ points from $S_r$, the remaining subgraph contains at least one edge. Hence, its chromatic number is $2$. Therefore,
    
    \[\Br(1,2)\geq \Br_{S_r}(1)>r.\]
    
    Since $r$ was arbitrary, it follows that $\Br(2,1)$ is not finite. Thus, by Proposition \ref{Prop:dim}, $\Br(1,d)$ is not finite for every $d>1$. Therefore, $\Br(1,d)=1$ if and only if $d=1$.
\end{proof}

\begin{proposition}\label{Prop:Borsuk}
     Let $d$ be a fixed positive integer. Then $\Br(k, d)=0$ for all $k\geq \beta_d$, where $\beta_d$ is the Borsuk number in dimension $d$.
\end{proposition}

\begin{proof}
    Since for every finite point set $S\subset \mathbb{R}^d$, $\chi(\Bor(S))=\beta(S)\leq \beta_d$, then $\Br_S(\beta_d)=0$. Hence
    $\Br(\beta_d, d)=0$. By Proposition \ref{Prop:color}, $\Br(k, d)=0$ for every $k\geq\beta_d$.
\end{proof}

\begin{remark}
    In dimension $2$, Borsuk proved in \cite{borsuk1933drei} that the Borsuk number is $3$. Hence, by Proposition~\ref{Prop:Borsuk}, we have $\Br(k,2)=0$ for all $k\geq 3$. Later, Dol’nikov showed in Theorem 1 \cite{dol2000some} that at most one vertex must be removed  from an arbitrary diameter graph to obtain a bipartite subgraph. Consequently, by Proposition \ref{Prop-simplex}, $\Br(2,2)=1.$
\end{remark}
\bigskip
\begin{remark}
    In dimension $3$, we have that $\Br(2,3)$ is not finite by Alon's construction in \cite{alon2003problems} as we mentioned in the introduction. Theorem \ref{Thm:Brass} and Proposition \ref{Prop-simplex} establishes that $\Br(3,3)=1$. Finally, for $k\geq 4$, $\Br(k,3)=0$ by Propositions \ref{Prop:color} and \ref{Prop:Borsuk}.
\end{remark}

\begin{comment}
    
\begin{theorem} [Theorem 8.1 \cite{alon2003problems}, 2003]
    For every positive integer $s$ there is a finite ser of points $X=X(s)$ in $\mathbb{R}^3$ that cannot be partitioned into three sets $A, B$ and $C$ such that $A$ is of size at most $s$ and the diameter of each of the sets $B$ and $C$ is strictly smaller than that of $X$. 
\end{theorem}
\end{comment}

\subsection{Higher dimensions}

Since the Borsuk number $\beta_4$ has not been exactly computed yet, it is not clear from which $k$ the Brass number $\Br(k,4)$ is $0$. Recently, Tolmachev and Voronov claim in \cite{tolmachev2026reducingupperboundborsuk}
that $\beta_4\leq8$, which would imply that $\Br(k,4)=0$ for every $k\geq8$. 

On the other hand, we know that $\Br(1,4)$ and $\Br(2,4)$ are not finite by Propositions \ref{Prop:dim} and \ref{Prop:1-color}. Then, the unknown Brass numbers in dimension 4 would be for $k=3,\dots, 7$.

We believe that $\Br (3,4)$ is not finite, although we have not found a construction that proves this fact. However, we consider the number $\Br (4,4)$ a more interesting case, since its understanding might give us a better idea of whether the Borsuk number $\beta_4$ is equal to 5 (as the Borsuk conjecture states) or not. In case $\Br (4,4)$ is not finite, we believe that it might be possible to extend such a construction to create a counterexample for the Borsuk conjecture in dimension 4.

\section{Final comments} \label{Sec:Final}

On our quest to prove that $\Br (3,4)$ is not finite, and after several attempts to find such family of graphs, we identified a pattern that we state in the following question.

\bigskip

\begin{question}\label{Question:criticalR4}
Let $S_1, S_2\subset \mathbb{R}^d$ be unit diameter finite point sets, such that their union has unit diameter as well. If both $\Bor{(S_1)}$ and $\Bor{(S_2)}$ are $d$-critical graphs, is true that $|S_1\cap S_2|\geq d-2$?.
\end{question}

Note that when $\Bor{(S_1)}$ and $\Bor{(S_2)}$ are $d$-cliques in $\mathbb{R}^d$, Question \ref{Question:criticalR4} has a positive answer due Theorem 1 in. Also notice that for $d=3$, Question \ref{Question:criticalR4} has a positive answer due Theorem 2 in \cite{dol2000some}. 

Finally, we would like to conclude this paper by adding a table that summarizes the known and unknown Brass numbers up to dimension 5.
\begin{table}[h!]
    \centering
    \begin{tabular}{|c|c|c|c|c|c|}
    \hline
         $k\backslash d$ & $1$ &$2$ & $3$& $4$ &$5$ \\\hline
         1& 1  & $\infty$ &$\infty$   &$\infty$   &$\infty$\\
         2& 0& 1 &$\infty$ & $\infty$& $\infty$ \\
         3& 0 &0 & 1& $\infty$& $\infty$\\
         4 & 0 & 0 & 0 & unknown & unknown\\
         5 & 0 & 0 & 0 & unknown & unknown \\
         6 & 0 & 0 & 0 & unknown & unknown \\
         7 & 0 & 0 & 0 & unknown & unknown \\
         8 & 0 & 0 & 0 & 0 & unknown \\
         \hline
    \end{tabular}
    \caption{Brass numbers $\Br (k,d)$.}
    \label{tab:placeholder}
\end{table}

\section{Declaration of AI usage}
\begin{table}[h!]
    \centering
    \begin{tabular}{c c}
    \hline
        Activity & AI used \\
        \hline
         Asking the research questions & No\\
         Coming up with the approach & No\\
         Proof development & No\\
         Writing and exposition & Yes\\
         Checking for mistakes & Yes\\
         Images & Yes\\
         \hline
    \end{tabular}
    \caption{Activities developed by AI}
    \label{tab:AI}
\end{table}

The constructions, proofs and ideas presented in this paper are due to the author. The main contribution of AI in this paper was in the proof of Lemma \ref{lemma:alpha}, where AI helped the author to simplify the proof, suggesting the cylinder grid to have a clearer counting.

\bibliography{main}
\end{document}